\documentclass{amsart}
\usepackage[english]{babel}
\usepackage{amssymb} 
   \newtheorem{Lem}{Lemma}[section]
    \newtheorem{theorem}[Lem]{Theorem}
    \newtheorem{conjecture}[Lem]{Conjecture}
    \newtheorem{proposition}[Lem]{Proposition}

   \theoremstyle{definition}
    
    \newtheorem{question}[Lem]{Question}
    \newtheorem{example}[Lem]{Example}
    \newtheorem{remark}[Lem]{Remark}

\renewcommand{\AA}{\mathbb A}
\newcommand{\QQ}{\mathbb Q}
\newcommand{\CC}{\mathbb C}
\newcommand{\FF}{\mathbb F}
\newcommand{\PP}{\mathbb P}
\newcommand{\ZZ}{\mathbb Z}
\newcommand{\HH}{\mathbb H}
\newcommand{\TT}{\mathbb T}
\newcommand{\RR}{\mathbb R}
\newcommand{\cO}{\mathcal O}
\newcommand{\cS}{\mathcal S}
\newcommand{\cM}{\mathcal M}
\newcommand{\cW}{\mathcal W}

\DeclareMathOperator{\Spec}{Spec}
\DeclareMathOperator{\etal}{\acute{e}t}
\DeclareMathOperator{\trace}{tr}
\DeclareMathOperator{\Fr}{Fr}
\DeclareMathOperator{\Fix}{Fix}

\DeclareMathOperator{\Gal}{Gal}
\DeclareMathOperator{\GL}{GL}
\DeclareMathOperator{\SL}{SL}
\DeclareMathOperator{\End}{End}
\DeclareMathOperator{\Sym}{Sym}

\DeclareMathOperator{\Art}{Art}
\DeclareMathOperator{\ab}{ab}
\DeclareMathOperator{\Kthree}{K3}
\DeclareMathOperator{\NS}{NS}
\DeclareMathOperator{\Km}{Km}
\DeclareMathOperator{\Nm}{Nm}
\DeclareMathOperator{\diag}{diag}

\begin{document}

\title{Modularity of Calabi-Yau varieties}
\author{Klaus Hulek \and Remke Kloosterman \and Matthias Sch\"utt}
\address{
Institut f\"ur Algebraische Geometrie,
Universit\"at Hannover,
Welfengarten 1, 30167 Hannover, Germany}
\email{hulek@math.uni-hannover.de}

\address{Institut f\"ur Algebraische Geometrie,
Universit\"at Hannover,
Welfengarten 1, 30167 Hannover, Germany}
\email{kloosterman@math.uni-hannover.de}

\address{
Institut f\"ur Algebraische Geometrie,
Universit\"at Hannover,
Welfengarten 1, 30167 Hannover, Germany,}
\email{schuett@math.uni-hannover.de}

\begin{abstract} In this paper we discuss recent progress on the modularity of Calabi-Yau varieties. We focus mostly on the case of surfaces and threefolds. We will also discuss some progress on the structure of the $L$-function in connection with mirror symmetry. Finally, we address some questions and open problems.
\keywords{Modular Calabi-Yau varieties, Mirror symmetry}
\end{abstract}

\thanks{
We cordially thank the DFG who has generously supported our research in
the framework of the Schwerpunktprogramm 1094 `Globale Methoden in der
komplexen Geometrie'. In our project on the modularity of Calabi-Yau
varieties we closely cooperated with the research group in Mainz led by
D.~van~Straten. We thank him and C.~Meyer for numerous discussions. We
are grateful to S.~Cynk, J.~Top and H.~Verrill who are coauthors of some
of the papers discussed here and with whom we have collaborated in many
ways. Last but not least we thank S.~Kadir and N.~Yui for many stimulating
discussions.
}
\maketitle

\section{Introduction}\label{intro}
In this paper we discuss recent progress on the modularity of Calabi-Yau varieties. We first explain what we mean by modularity. A $d$-dimensional variety $X/\QQ$ is called modular if the $L$-function of the Galois-representation on $H^d_{\etal}(X,\QQ)$ equals the product of $L$-functions of modular forms up to the factors associated to the bad primes. In practice, if $d$ is even then $H^d_{\etal}(X,\QQ)$ contains many classes of $d/2$-dimensional cycles. The Galois representation on the subspace generated by these classes is very often easy to calculate. In this case we consider only the subrepresentation on the orthogonal complement of the images of the algebraic cycles.

Almost trivially, the modular curves associated to subgroups of $\SL(2,\ZZ)$ of finite index are modular. The first non-trivial example are the elliptic curves over $\QQ$ with complex multiplication. Deuring \cite{Deu} proved that their $L$-series is the $L$-series of a Gr\"ossencharacter.
One can associate a modular form to such an $L$-series by the work of Hecke.

It has been a long standing conjecture that all elliptic  curve over $\QQ$ are modular. This conjecture was known as the Taniyama-Shimura-Weil conjecture. Wiles \cite{Wi} proved a large part of this conjecture. A few years later Breuil, Conrad, Diamond and Taylor \cite{BCDT} completed the proof.
In this paper we will concentrate on  Calabi-Yau varieties of dimension two and higher. 

In  Section~\ref{section1} we discuss some preliminaries, such as the Weil conjectures, \'etale cohomology and the definition of an $L$-function of a Galois representation. In Section~\ref{sec_Auto} we explain how one can associate a Galois representation to a modular form. Furthermore, we discuss  several general conjectures concerning the modularity and automorphicity of Galois representations, 
namely the Fontaine-Mazur conjectures and we very briefly describe the Langlands' programme.

In Section~\ref{sec_2dimensional} we come back to the elliptic curve case in more detail and discuss a similar result for rigid Calabi-Yau threefolds. In Section~\ref{sec_examples} we discuss some examples of modular varieties, where we mainly focus on singular $K3$-surfaces and Calabi-Yau threefolds, although we also consider some higher dimensional examples.

In Section~\ref{sec_Mirror} we discuss some recent progress on determining the structure of the $L$-function of  Calabi-Yau varieties and the $L$-function of the mirror partner. Here we sketch a new approach to compute the
zeta function of the one-parameter quintic family and its mirror.
Finally, in Section~\ref{Sec_open} we address some open questions.
\section{Varieties over $\QQ$.}\label{section1}

We shall discuss first varieties which are defined over the field $\QQ$ of rational numbers.
Later we shall also, more generally, consider varieties defined over number fields.
Let $X$ be a $d$-dimensional projective variety defined over $\QQ$, i.e., a variety defined by the vanishing of a finite 
number of homogeneous
polynomials with rational coefficients. 

Let $p$ be a prime number. We say that $X$ has \emph{good reduction} at $p$ if there exists a variety $X'/\QQ_p$ such that $X$ and $X'$ are isomorphic over $\QQ_p$ and the reduction of $X'$ modulo $p$ is smooth. If  this is the case then we take $X_p$ to be the reduction of $X'$. We call such $p$ a \emph{good prime} (for $X$). If $p$ does not satisfy the above mentioned property then we say that $X$ has \emph{bad reduction} at $p$. Such a $p$ is called a \emph{bad prime}. The number of bad primes is always finite.

If we fix a model of $X$ over $\ZZ$ then for all but finitely many primes $p$ we have a another way to define $X_p$, i.e., we have that $X_p$ is isomorphic to  $X\times_{\Spec \ZZ} \FF_p$, provided that the latter  is smooth.

If $E$ is an elliptic curve over $\QQ$ we can find a model over $\ZZ$ such that for all good primes $p$ the curve $E\times_{\Spec \ZZ} \FF_p$ is smooth. However, in general one cannot expect to find one model that suffices for all good primes.

If we consider a variety $X$ defined over a field $K$ that is not algebraically closed, then we indicate by $\overline{X}$ the same variety considered over the algebraic closure $\overline{K}$.
\subsection{The zeta function}\label{sub_zeta}
Let $q=p^m$ be a prime power and let $X/\FF_q$ be  a variety. Set 
$$
N_{q^r} = \# X(\FF_{q^r})
$$
be the number of points of $X$ over $\FF_{q^r}$. Then the zeta function of $X/\FF_q$ is defined by
$$
Z_q(t)= \exp\left( \sum_{r=1}^{\infty} N_{q^r} \frac{t^r}{r}\right) \in \QQ[[t]].
$$
Weil \cite{WeilC} has made a series of famous conjectures concerning this function:

\begin{theorem}[Weil conjectures] 
Let $X$ be a smooth $d$-dimensional projective variety defined over the field $\FF_{q}$, $q=p^r$. 
The function $Z_q(t)$ satisfies the following properties
\begin{enumerate}
\item \emph{Rationality}, i.e.
$$
Z_q(t)= \frac{P_q(t)}{Q_q(t)} \mbox{ for polynomials } P_q(t), Q_q(t) \in \ZZ[t]. 
$$
\item The function $Z_q(t)$ satisfies a \emph{functional equation}
$$
Z_q\left(\frac{1}{q^dt}\right)= \pm q^{de/2} t^{e}Z_q(t)
$$
where $e$ is the self-intersection of the diagonal in $X \times X$.
\item The \emph{Riemann hypothesis} holds, i.e.
$$
Z_q(t)= \frac{P_{1,q}(t)\cdots P_{2d-1,q}(t)}{P_{0,q}(t)\cdots P_{2d,q}(t)}
$$
where $P_{0,q}(t)=1-t, P_{2d,q}(t)=1-q^dt$ and
$$
P_{i,q}(t)= \prod_{j=1}^{b_i}(1- \alpha_{ij}t)
$$
for $1 \leq i \leq 2d-1$, and where the $\alpha_{ij}$ are algebraic integers of complex absolute 
value $\mid \alpha_{ij}\mid=q^{i/2}$.
\item Assume that $X$ arises as the reduction of a variety
defined over a number field $K$. Then the $b_i$ are the topological \emph{Betti numbers} of the complex
variety ${X}_{\CC}$ and $e = \sum_{i=0}^{2d}(-1)^ib_i$ is the Euler number.
\end{enumerate}
\end{theorem}

The most difficult part of these conjectures is the Riemann hypothesis (in Subsection~\ref{sub_Lfunctions} we explain
why this is related to the classical Riemann hypothesis), 
which was finally proved by Deligne \cite{DelRH} in
1974. These conjectures, especially point 4, establish a relationship between the arithmetic and the
topology of an algebraic variety.

\subsection{\'Etale cohomology and Frobenius}\label{sub_Frobenius}
Deligne proved the Weil conjectures in the technical framework of \'etale cohomology.

This cohomology theory was developed by Grothendieck with a view towards proving the Weil
conjectures.
We refer the reader to \cite{FK} for an account of this theory, here we simply state its basic properties.
Let $Y$ be a smooth projective variety defined over an algebraically closed field
of characteristic $p$, which can be either $0$ or positive (we shall use this later for 
$Y=\overline{X}_p$ and $Y=\overline{X}$). Choose some prime $l \neq p$. 
Then there is a cohomology theory, called  
\emph{\'etale cohomology}, which associates to $Y$ certain $\QQ_l$-vector spaces
$H^i_{\etal}(Y,\QQ_l)$ for $i \geq 0$, which have many properties  similar to the 
classical (singular) 
cohomology in characteristic $0$. In particular, the \'etale cohomology groups have the following properties:
\begin{enumerate}
\item If $Y$ has dimension $d$, then $H^i_{\etal}(Y,\QQ_l)=0$ for
$i > 2d$. 
\item All $\QQ_l$-vector spaces $H^i_{\etal}(Y,\QQ_l)$ 
are \emph{finite dimensional}.
\item The cohomology groups behave \emph{functorially} with respect to morphisms of 
smooth algebraic
varieties.
\item \emph{Poincar\'e duality} holds, i.e. there is an isomorphism
$H^{2d}_{\etal}(Y,\QQ_l) \cong \QQ_l$ and for each $i \leq d$ a
perfect  pairing
$$
H^i_{\etal}(Y,\QQ_l) \times H^{2d-i}_{\etal}(Y,\QQ_l) \to
H^{2d}_{\etal}(Y,\QQ_l) \cong \QQ_l.
$$
\item The \emph{K\"unneth decomposition} theorem for products holds.
\item The \emph{Lefschetz fixed point formula} holds: Let $f: Y \to Y$ be a morphism such that the set  of fixed points $\Fix(f)$ is finite and  $1-\mathrm{d}f$ is injective (here $\mathrm{d}f$ denotes 
the differential of $f$),
then
$$
\# \Fix(f)= \sum_{i=0}^{2d}(-1)^i \trace(f^*\mid H^i_{\etal}(Y,\QQ_l)) .
$$  
\item A \emph{comparison theorem} with singular 
cohomology holds: if $Y$ is smooth and projective over $\CC$,
then there are isomorphisms of $\CC$-vector spaces
$$
H^i_{\etal}(Y,\QQ_l) \otimes_{\QQ_l} \CC \cong  H^i(Y,\CC).
$$
\end{enumerate}

The crucial property for our purposes is the Lefschetz fixed point formula. The condition that
$1 -\mathrm{d}f$ is injective means that the fixed points have multiplicity one.
Assume now that
$X$ is a smooth projective variety defined over $\QQ$ and that $p$ is a prime of good reduction.
Consider the geometric Frobenius morphism $\Fr_p: \overline{X}_{p} \to \overline{X}_{p}$ given by 
sending the coordinates to their $p$-th power. Then 
$$
X(\FF_{p^r})=\Fix(\Fr_p^r).
$$
Now choose a different prime $l \neq p$ (the choice of this prime will not matter). 
Then by the Lefschetz fixed point formula
$$
N_{p^r} = \#X(\FF_{p^r})=\sum_{i=0}^{2d}(-1)^i 
\trace((\Fr_p^r)^*\mid H^i_{\etal}(\overline{X}_{p},\QQ_l)) .
$$
We define polynomials
\begin{equation}\label{def_P}
P_{i,p^r}(t)= \det(1-(\Fr_{p^r})^*t\mid H^i_{\etal}(\overline{X}_{p},\QQ_l)).
\end{equation}
Using the identity 
$$
\sum_{n \geq 1} \trace(f^n\mid V) \frac{t^n}{n}= - \log \det(1-ft) 
$$ 
for an endomorphism $f$ of a finite-dimensional vector space $V$, 
it is then essentially an exercise in formal power series to show that 
$$
Z_{p^r}(t)= \frac{P_{1,p^r}(t)\ldots P_{2d-1,p^r}(t)}{P_{0,p^r}(t)\ldots P_{2d,p^r}(t)}.
$$
In particular, the polynomials $P_{i,p^r}(t)$ introduced in Section~\ref{sub_zeta} coincide
with the polynomials from equation (\ref{def_P}).

Deligne's proof of the Riemann hypothesis for function fields of varieties defined over
finite fields then follows from
\begin{theorem}[{Deligne, \cite{DelRH}, \cite{DelRHa}}] Let $Y$ be a smooth projective variety defined over the
field $\FF_q$, where $q=p^r$ is a prime power. 
For all $i$ and $l\neq p$ the eigenvalues of $\Fr_q^*=(\Fr_p^r)^*$ on
$H^i_{\etal}(\overline{Y}_p,\QQ_l)$ are algebraic integers of absolute 
value $q^{i/2}$.
\end{theorem}

In the case of a curve $X$ there is another, more explicit, construction of the \'etale cohomology. Let $J(X)$ denote the Jacobian of $X$, then
\[ H^1_{\etal}(X,\QQ_l) \cong \left(\lim_{\stackrel{\longleftarrow}{n}}  J(X)[l^n]\right) \otimes_{\ZZ_l} \QQ_l.\]
If $X$ is defined over $\QQ$ then there is a natural action of $\Gal(\overline{\QQ}/\QQ)$ on $J(\overline{X})[l^n]$ which yields a Galois action on $H^1_{\etal}(\overline{X},\QQ_l)$. For a general variety $X$ defined over $\QQ$ there is a natural action of $\Gal(\overline{\QQ}/\QQ)$ on $H^i_{\etal}(\overline{X},\QQ)$, but this nice description is missing.

\subsection{$L$-functions and Galois representations}\label{sub_Lfunctions}
We are now in a position to discuss $L$-functions. Let $X$ be a model of a variety defined over $\QQ$. Whereas we have, so far, considered the
reduction of a variety $X$ defined over $\QQ$ for a fixed prime $p$, we shall now vary 
the prime $p$. For any $i \geq 0$ we define the $L$-function
$$
L(H^i_{\etal}(\overline{X},\QQ_l),s)= (*) \prod_{p\in \mathcal{P}}\frac{1}{P_{i,p}(p^{-s})},
$$ 
where $\mathcal{P}$ is the set of primes of good reduction and $(*)$ denotes suitable Euler factors for the 
primes of bad reduction. In this survey article we shall not discuss the (difficult) question of
how to define the
Euler factors for bad primes in general, though we will explain this for the case of
elliptic curves (see Section~\ref{sub-Wiles}). Note that if $X=\PP^1$ and $i=0$, then we obtain exactly  
Riemann's zeta function $\zeta(s)$, and if $X$ is a smooth projective curve, one thus obtains the 
Hasse zeta function
of $X$. 

In most cases, the interesting cohomology of a variety $X$ is the middle cohomology.
For this reason we define the $L$-series of a $d$-dimensional variety $X$ as 
$$
L(X,s):= L(H^d_{\etal}(\overline{X},\QQ_l),s).
$$
 
We have already remarked that 
$H^i_{\etal}(\overline{X},\QQ_l) \otimes_{\QQ_l} \CC \cong H^i(X_{\CC},\CC)$
as $\CC$-vector spaces. 

Let $p$ be a good prime. A further comparison theorem says that the reduction map gives an isomorphism
\[ \psi: H^i_{\etal}(\overline{X},\QQ_l) \stackrel{\sim}{\longrightarrow}
H^i_{\etal}(\overline{X}_p,\QQ_l)\] 
as $\QQ_l$-vector spaces.
The isomorphism $\psi$ is also compatible with the  Galois module structure of these spaces in the following sense.
The Galois group 
$\Gal(\overline{\FF}_p/\FF_p)$ acts on $\overline{X}_p$ thus making
$H^i_{\etal}(\overline{X}_p,\QQ_l)$ a $\Gal(\overline{\FF}_p/\FF_p)$-module.
Similarly, the absolute Galois group $\Gal(\overline{\QQ}/\QQ)$ acts on  
$\overline{X}$ giving $H^i_{\etal}(\overline{X},\QQ_l)$ a 
$\Gal(\overline{\QQ}/\QQ)$-module structure. This makes Poincar\'e duality a Galois-equivariant pairing.

For every embedding 
$\overline{\QQ}\rightarrow \overline{\QQ_p}$ one obtains an embedding \[\Gal({\overline{\QQ_p}}/\QQ_p)\hookrightarrow \Gal(\overline{\QQ}/\QQ).\] Two different embeddings yield conjugated subgroups of $\Gal(\overline{\QQ}/\QQ))$. Let $\widetilde{\Fr}_p$ be a lift of $\Fr_p$. Then the action of
$\psi \circ \widetilde{\Fr}_p$ on $H^i_{\etal}(\overline{X}_p,\QQ_l)$ is conjugated to the action of $\Fr_p$.

Let
$$
\rho_i(X): \Gal(\overline{\QQ}/\QQ) \to \GL(H^i_{\etal}(\overline{X},\QQ_l)).
$$
be the natural Galois representations.
Since any two lifts of $\Fr_p$ are conjugate
they have the same trace and determinant. In this way we can still speak about the trace and determinant of 
Frobenius at a good prime $p$ in connection with $\Gal(\overline{\QQ}/\QQ)$-representations.

\section{Automorphic Origin of $l$-adic Representations}\label{sec_Auto}

\subsection{The Galois representation associated to a modular form} \label{subsec_Galoisassmodform}
To recall the notion of (elliptic) modular form, we recall first that the group
$\SL(2,\ZZ)$ acts on the complex upper half plane $\HH$ by
$$
\begin{pmatrix} 
a & b\\ c&d 
\end{pmatrix}: \tau \mapsto \frac{a \tau + b}{c \tau + d}.
$$
We can define $\HH^*$ by adding to $\HH$ a copy of $\PP^1(\QQ)$ in the following way: the point $(x:1)$ is 
identified with $x\in \QQ\subset \CC$ and the point $(1:0)$ is identified with the point at 
infinity along the complex axis, usually denoted by $i\infty$. 
The action  of $\SL(2,\ZZ)$ extends to $\HH^*$.

For a given integer $N$ we define the group
$$
\Gamma_0(N)= 
\left\{ \begin{pmatrix} 
a & b\\ c&d 
\end{pmatrix} \in \SL(2,\ZZ); \, c \equiv 0 \bmod N \right\}.
$$
The quotient $Y_0(N) = \Gamma_0(N) \backslash \HH$ is a Riemann surface, which can be compactified to
a projective curve $X_0(N)=\Gamma_0(N)\backslash \HH^*$. The curve $X_0(N)$ is called the \emph{modular curve 
of level $N$} and
the points in $X_0(N)\setminus Y_0(N)$ are called \emph{cusps}. 

An elliptic modular form of \emph{weight} $k$ and \emph{level} $N$ is a holomorphic
function
$$
f: \HH \to \CC
$$
which satisfies the following properties:
\begin{enumerate}
\item For matrices $ M=\begin{pmatrix} a & b\\ c&d 
\end{pmatrix} \in \Gamma_0(N)$ the function $f$ transforms as
$$
f\left(\frac{a \tau + b}{c \tau + d}\right) = (c \tau + d)^k f(\tau).
$$
\item The function $f$ is holomorphic at the cusps.
\end{enumerate}
To explain the latter condition, we consider one of the cusps, namely the cusp given by 
adding $i \infty$. 
The matrix 
$ \begin{pmatrix} 1 & 1\\ 0&1  \end{pmatrix}$ is an element of $\Gamma_0(N)$ and we thus have 
$$
f(\tau+1)= f(\tau).
$$
Writing $q = \exp^{2\pi i \tau}$ we can, therefore, consider the Fourier expansion
$$
f(\tau)=f(q) = \sum_n b_nq^n.
$$
The cusp $i\infty$ corresponds to $q=0$ and $f$ is holomorphic at the cusp $i \infty$ if $b_n=0$ for $n<0$. We say that $f$ vanishes at this cusp, if in addition
$b_0=0$. A \emph{cusp form} is a modular form which vanishes at all cusps. 
We shall also consider modular forms with a $\bmod N$-Dirichlet character $\chi$. In this case the transformation behaviour
given in the definition of modular forms is
$$
f\left(\frac{a \tau + b}{c \tau + d}\right) = \chi(d)(c \tau + d)^k f(\tau).
$$
The vector space of weight $k$ cusp forms for $\Gamma_0(N)$ is denoted by $S_k(\Gamma_0(N))$ and,
similarly, the vector space of weight $k$ cusp forms with character $\chi$ is denoted 
by $S_k(\Gamma_0(N),\chi)$. 

On $S_k(\Gamma_0(N),\chi)$ we can define operators $T_p$ for prime numbers $p\nmid N$. For a precise definition  we refer to \cite[page 244]{Kn}.
The operators $T_p$ are called \emph{Hecke operators}. 
They generate a subalgebra of $\End(S_k(\Gamma_0(N),\chi)$ called the \emph{Hecke algebra} $\TT$. 

Since Hecke operators $T_q$ and $T_q$ commute  for distinct primes $p,q$ not dividing $N$, 
we have simultaneous eigenspaces for all $T_p$. 
A form that is an eigenvector for all $T_p$ is called a \emph{Hecke eigenform}. Let $g\in S_k(\Gamma_0(N),\chi)$ be a Hecke eigenform, so $T_p(g)=c_pg$. Then there is a minimal $N' \mid N$ and a unique eigenform  
\[ f=\sum b_n q^n\in S_k(\Gamma_0(N'),\chi)\]
with the same system of eigenvalues such that $b_p=c_p$ for all $p\nmid N$. The eigenform $f$ is called a (normalised) newform.

Let $f$ be a normalised Hecke eigenform of some level $N$ and weight $k>1$. 
We will now explain how one can associate to $f$  a Galois representation 
\[\rho_f:\Gal(\overline{\QQ}/\QQ)\rightarrow \GL(2,\QQ_{l})\]
in the case that  $k>2$. The case $k=2$ is treated in Remark~\ref{rem_repgeo}. The following construction is due to Delinge \cite{DelCon}.
Using the trivial inclusion $S_k(\Gamma_0(N),\chi)\subset S_k(\Gamma_0(NM,\chi)$ we may assume that $N>4$.
 
In order to define $\rho_f$ we start by considering the
subgroup $\Gamma_1(N)\subset \Gamma_0(N)$ of matrices $M$ such that 
\[M\equiv \left( \begin{array}{cc} 1 &\overline{b} \\ 0&1\end{array} \right)\bmod N.\]
A modular form in $S_k(\Gamma_0(N),\chi)$ is also a modular form for the group $\Gamma_1(N)$. 
Note, however, that
the construction we present below uses the group  $\Gamma_1(N)$ in an essential way.
As above, we can define modular curves $Y_1(N)$ and $X_1(N)$. 
One can show that both curves have a natural model over $\QQ$.
Since we assumed that $N>4$ we have that $Y_1(N)$ is a {\em fine} moduli-space for pairs $(E,P)$ with $E$ 
an elliptic curve 
and $P$ a point of order $N$. Let $\pi: \mathcal{E}_1(N) \rightarrow Y_1(N)$ be the associated universal family.
We define define the $l$-adic 
sheaf $\mathcal{F}_k=\Sym^{k-2} R^1\pi_*\QQ_l$ on the \'etale site of $Y_1(N)$.

Let $j: Y_1(N)\rightarrow X_1(N)$ be the natural morphism. Then the \'etale cohomology group
\[H:=H^1_{\etal}(\overline{X_1(N)},j_*\mathcal{F}_k)^{\vee{}}\] 
has an action of the Hecke-algebra $\TT\otimes_{\QQ}\QQ_l$. It turns out that $H$ is a direct sum of 
free rank two Hecke-algebras.  Recall that we had fixed a normalised Hecke eigenform $f$. Let $b_p$ denote the eigenvalue of $T_p$. Consider the morphism  $\TT \otimes_{\QQ}\QQ_l\rightarrow\overline{\QQ_l}$ 
defined by $T_p \rightarrow b_p$. Define
\[ V_f := H^1_{\etal}(\overline{X_1(N)},j_*\mathcal{F}_k)^{\vee{}} 
\otimes_{\TT\otimes_{\QQ} \QQ_l} \overline{\QQ_l}.\]
It turns out that $V_f$ is two-dimensional and that there is a natural 
action of $\Gal(\overline{\QQ}/\QQ)$ on $V_f$. This gives our representation $\rho_f$.

\begin{remark}\label{rem_repgeo} If $f$ has weight $k=2$ then $f\frac{dq}{q}$ induces a holomorphic one-form on $X_1(N)$. 
Consider the two-dimensional subspace 
\[ V'_f:=f\frac{dq}{q}\CC\oplus \overline{f \frac{dq}{q}}\CC\subset H^1(\overline{X_1(N)},\CC)\]
Using the comparison theory of \'etale cohomology and singular cohomology we obtain a two-dimensional 
subspace $V_f$ of $H^1_{\etal}(\overline{X_1(N)},\overline{\QQ_l})$. It turns out that $V_f$ is Galois-invariant. 
One can show that this gives a definition of $\rho_f$.

Similarly, if $f$ has weight $k>2$ then $f\frac{dq_1}{q_1}\wedge \dots \wedge \frac{dq_{k-2}}{q_{k-2}}$ induces a holomorphic $k-2$-form on $\Sym^{k-2} \overline{\mathcal{E}_1(N)}$. This leads to a two-dimension subrepresentation of $H^{k-2}(\overline{\mathcal{E}_1(N)},\overline{\QQ_l})$.
This provides an alternative definition of $\rho_f$. Although this construction seems to be known to many experts, we could only find references for  the special cases $k=3$ \cite{ShiMod} and $k=4$ \cite{SY}.
\end{remark}

To a normalised Hecke newform $f$ we can associate an $L$-function $L(f,s)$. 
One way of defining  this is $L(f,s):=L(\rho_f,s)$, i.e., as the $L$-series of the Galois-representation $\rho_f$. Equivalently, we can do the following:
consider the Fourier expansion $f=\sum b_nq^n$. One can associate to this its
\emph{Mellin transform}
$$
L(f,s)= \sum_n b_nn^{-s}.
$$ 
If $f$ is a normalised Hecke newform with respect to the group 
$\Gamma_0(N)$, then the Fourier coefficients satisfy the properties
\begin{eqnarray*}
b_{p^r}b_p&=&b_{p^{r+1}}+p^{k-1}b_{p^{r-1}}  \mbox{ for } p  \mbox{ prime, } p \nmid N\\
b_{p^r}&=&(b_p)^r  \mbox{for } p  \mbox{ prime, } p\mid N \\
b_n b_m &=& b_{nm} \mbox{ if } (n,m)=1
\end{eqnarray*}
where $k$ is the weight of the form $f$.
It follows from this that the series $L(f,s)$ has  a product expansion
\begin{equation}\label{eqn_Lseriesmodform}
L(f,s) = \sum_{n \geq 1} {b_n}{n^{-s}}= \prod_p 
\frac{1}{1-b_pp^{-s} + \chi(p)p^{k-1-2s}}
\end{equation}
where $\chi(p)=0$ if $p \mid N$.
\subsection{Fontaine-Mazur conjecture}\label{sub_FM}
Here we shall briefly discuss the Fontaine-Mazur conjecture.
A more detailed discussion of this conjecture as well as many references for the facts mentioned below, 
can be found in Taylor's survey paper \cite{Tay}.

Consider a variety $X/\QQ$. We have already mentioned that the \'etale cohomology 
groups $H^i_{\etal}(\overline{X},\overline{\QQ_l})$ have a natural Galois action, i.e., 
there is a continuous homomorphism
\[ \rho: \Gal(\overline{\QQ}/\QQ)\rightarrow \GL(H^i_{\etal}(\overline{X},\overline{\QQ_l})).\]

In order to discuss Galois-representations we first have to list some subgroups of $\Gal(\overline{\QQ}/\QQ)$.
For each embedding $\overline{\QQ}\hookrightarrow \overline{\QQ_p}$ we obtain $\Gal(\overline{\QQ_p}/\QQ_p)$ as 
a subgroup of $\Gal(\overline{\QQ}/\QQ)$. Two such embeddings 
$\overline{\QQ}\hookrightarrow \overline{\QQ_p}$ give subgroups  of $\Gal(\overline{\QQ}/\QQ)$ that differ only up to conjugacy. Fix for every prime $p$ an embedding $\overline{\QQ}\hookrightarrow \overline{\QQ_p}$.

The group $\Gal(\overline{\QQ_p}/\QQ_p)$ maps onto $\Gal(\overline{\FF_p}/\FF_p)$. 
The kernel of this map is called the inertia subgroup at $p$, denoted by $I_p$.

Let $V$ be a finite dimensional vector space over $\QQ_l$. Let $\rho': \Gal(\overline{\QQ}/\QQ)\rightarrow\GL(V)$ be an  irreducible subrepresentation of $\rho$. Then the representation $\rho'$ satisfies the following properties:
\begin{enumerate}
\item $\rho'$ is only ramified at finitely many primes. I.e., $\rho'(I_p)=\{1\}$ for all but finitely many 
primes $p$.
\item The restriction of $\rho'$ to the subgroup $\Gal(\overline{\QQ_l}/\QQ_l)$ is of de Rham 
type. (For a definition of this see \cite{FOa}, \cite{FOb}. Note that this is quite restrictive.).
\item For all but finitely many primes $p$ we have that all roots $\alpha$ of the characteristic 
polynomial of $\rho(\Fr_p)$ on $V$ satisfy $\mid \alpha\mid =p^{i/2}$, with $\mid \cdot \mid $ the complex absolute value.
\end{enumerate}

In general, an irreducible $l$-adic Galois-representation
\[ \rho: \Gal(\overline{\QQ}/\QQ)\rightarrow \GL(V)\]
with $V$ a finite dimensional $\overline{\QQ_l}$ vector space, is called {\em geometric} if it 
satisfies the above properties $1$ and $2$. 

The $i$-th {\em Tate twist} of a Galois-representation $\rho$ is defined as
the Galois representation obtain by tensoring $\rho$ with $ \mu_{l^\infty}^{\otimes i}$ if $i>0$ and tensoring $\rho$ with $\mu_{l^\infty}^{\vee{} \otimes -i}$ if $i<0$. Here $\mu_{l^n}$ denotes the $l^n$-th roots of unity and $\mu_{l^\infty}= \underleftarrow{\lim} \mu_{l^n}$. 

\begin{conjecture}[{Fontaine-Mazur, \cite[Conjecture 1]{FM}}]  Let $V$ be a finite dimensional $\overline{\QQ_l}$-vector space and 
let $\rho:\Gal(\overline{\QQ}/\QQ)\rightarrow \GL(V)$ be an irreducible 
geometric $l$-adic representation. Then there exists a smooth projective variety $X/\QQ$ and an integer $i$ such 
that $\rho$ is  a Tate twist of an irreducible subrepresentation of the standard Galois 
representation on $H^i(\overline{X},\overline{\QQ_\ell})$. In particular $\rho$ satisfies property 3.
Moreover, the $L$-function $L(\rho,s)$ has a meromorphic
 continuation to $\CC$. 
\end{conjecture}
A Tate twist of an irreducible subrepresentation of the standard Galois representation on $H^i(\overline{X},\overline{\QQ_l})$ is called a {\em Galois representation coming from geometry}.

In the two-dimensional case they make a stronger conjecture:
\begin{conjecture}[{Fontaine-Mazur, \cite[Conjecture 3c]{FM}}]  Let $V$ be a two-di\-men\-sio\-nal $\overline{\QQ_l}$-vector space. Let $\rho: \Gal(\overline{\QQ}/\QQ)\rightarrow \GL(V)$ be an irreducible geometric Galois representation that is not a Tate twist of a  finite representation. Then $\rho$ is isomorphic to the Tate twist of a Galois representation associated to a modular form.\end{conjecture}

Later on we will introduce rigid Calabi-Yau threefolds. A rigid Calabi-Yau threefold $X$ has $b_3(X)=2$. Suppose $X$ is defined over $\QQ$. Then we have a two-dimensional Galois representation. The above conjecture states that the Galois representation on $H^3_{\etal}(\overline{X},\overline{\QQ_l})$ is coming from a modular form.
This conjecture is also stated in \cite{SY}.

\subsection{Langlands' programme}
Here we shall only give a very brief outline. 
For a more extended introduction to this subject we refer to \cite{Tay}.
In the previous section we have seen a relation between Galois-representations, varieties over $\QQ$ and 
$L$-series. In this section we introduce a fourth class of objects, namely automorphic representations.

 Let $\hat{\ZZ}:=\underleftarrow{\lim} \ZZ / N\ZZ$. Denote $\AA^\infty = \hat{\ZZ}\otimes \QQ$ and 
$\AA=\AA^\infty \times \RR$. As a ring $\AA^\infty$ is the subring $\prod_p \QQ_p$ consisting of 
elements $(x_p)$ such that $x_p\in\ZZ_p$ for all but finitely many primes $p$. 
Note that the topology of $\AA^\infty$ is different from the subspace topology of $\prod \QQ_p$.

Let $W_{\QQ_p}$ be the Weil-group of $\QQ_p$, that is, all $\sigma \in \Gal(\overline{\QQ_p}/\QQ_p)$ such that their image in $\Gal(\overline{\FF_p}/\FF_p)$ is a power of $\Fr_p$.
This is a dense subgroup of $\Gal(\overline{\QQ_p}/\QQ_p)$. Local class field theory gives 
a `natural' isomorphism
\[ \Art_p : \QQ_p^*\stackrel{\sim}{\longrightarrow} W_{\QQ_p}^{\ab},\]
where $*^{\ab}$ means the abelianised  group. Furthermore, global class field theory gives an isomorphism
\[ \Art: \AA^*/\QQ^*\RR^*_{>0} \stackrel{\sim}{\longrightarrow} \Gal(\overline{\QQ}/\QQ)^{\ab}\]
which is compatible with $\Art_p$. (A large part of) class field theory can be described in 
terms of the map $\Art$. 
The Langlands' programme can be considered as an attempt to generalise this 
to $\GL(n,\QQ) \backslash \GL(n,\AA)$.

Fix $n\geq 1$. An automorphic form is a smooth function 
\[\GL(n,\QQ)\backslash\GL(n,\AA)\rightarrow \CC\]satisfying certain conditions: 
for example, the translates of $f$ under $\GL(n,\hat{\ZZ})\times O(n)$ (a maximal compact 
subgroup of $\GL(n,\AA)$) form a finite-dimensional vector space. Other conditions concern a generalisation of the
vanishing of $f$ at cusps and  limit the growth of $f$.
For reasons of space we will not give a complete definition, see \cite{Tay}. 
To an automorphic form one can associate an infinitesimal character $H$, where $H$ is a multiset of $n$ 
complex numbers (i.e., an element of $\Sym^n \CC$). Fix such a multiset $H$. Let $\mathcal{A}:=\mathcal{A}^o_H(\GL(n,\QQ)\backslash \GL(n,\AA))$ 
be the vector space of cuspidal automorphic forms with infinitesimal character $H$. We would like to 
consider $\mathcal{A}$ as a representation of $\GL(n,\AA)$. It turns out that $\GL(n,\RR)$ does not 
fix $\mathcal{A}$. However, $\mathcal{A}$ has an action of $\GL(n,\AA^\infty)\times O(n)$ 
and an action of 
the Lie-algebra $\mathfrak{gl}_n$. An irreducible constituent of $\mathcal{A}$ is called a 
{\em cuspidal automorphic representation} of $\GL(n,\AA)$. 

Consider $n=1$. One easily describes $\mathcal{A}^o_{\{s\}}(\QQ^*\backslash \AA^*)$, with $s\in \CC$,
in terms of 
$\mathcal{A}^o_{\{0\}}(\QQ^*\backslash\AA^*)$. It turns out that elements in this space are locally 
constant functions $f$ on
\[ \hat{\ZZ}^*\stackrel{\sim}{\longrightarrow} \AA^*/\QQ^* \RR^*_{>0}.\]
Thus $\mathcal{A}^o_{\{0\}}(\QQ^*\backslash \AA^*)$ is generated by continuous characters 
$\hat{\ZZ}^*\rightarrow \CC$.

Consider $n=2$. Write $H=\{s,t\}$. Then $\mathcal{A}\neq (0)$ only if $s-t\in i \RR$, $s-t\in \ZZ$ or 
$s-t\in (-1,1)$.
If $s-t\in \ZZ_{>0}$ then to every element in $\mathcal{A}$ we can associate a weight $1+s-t$ modular form and one can 
determine its level  in terms of $\mathcal{A}$. This shows that automorphic forms are 
generalisations of modular forms.

Suppose $n>1$. As in the case of modular forms we can associate local  $L$-series $L_p(\pi,s)$ to an 
automorphic form $\pi$. (Basically, we restrict the automorphic form $\pi$ to 
$\pi_p:\GL(n,\QQ_p)\subset \GL(n,\AA)$ in order to define a local $L$-factor $L_p(\pi_p,s)$.) 
The Euler product $\prod_p L_p$ converges on some half-plane, and has a meromorphic continuation to $\CC$.

A weak form of the Langlands' programme is to prove that every $L$-series of an automorphic form is also 
the $L$-series of an irreducible Galois-re\-pre\-sen\-tat\-ion coming from geometry and that this correspondence 
is `natural'.

\section{Two-dimensional Galois representations}\label{sec_2dimensional}

Wiles' proof of a special case of the Taniyama-Shimura-Weil conjecture was the essential ingredient in his proof of
Fermat's Last Theorem. In this section we briefly recall his result and discuss generalisations to
higher dimension, in particular to rigid Calabi-Yau threefolds.

\subsection{Elliptic curves and Wiles' theorem}\label{sub-Wiles}

We consider elliptic curves over $\QQ$, which we can assume to be in Weierstrass form
$$
y^2=x^3 + Ax + B,  \quad \Delta=4A^3+27B^2 \neq 0
$$
where $A,B \in \QQ$. For good primes $p$ the Lefschetz fixed point formula reads
$$
N_p = \#E(\FF_p) = 1 - \trace(\Fr_p^*\mid H^1_{\etal}(\overline{E}_{p},\QQ_l)) +p.
$$
Writing 
$$
a_p= \trace(\Fr_p^*\mid H^1_{\etal}(\overline{E}_{p},\QQ_l))
$$
this formula simply becomes
$$
N_p= 1 - a_p +p.
$$
The Taniyama-Shimura-Weil conjecture (TSW), in one of its forms, predicts that the 
sequence $N_p$, resp. $a_p$, is given (for good primes) by the
Fourier coefficients of 
a modular form. More precisely, 
\begin{theorem}[{Wiles-Taylor-Breuil-Conrad-Diamond, \cite{BCDT}}]\label{Thm_Wiles}
Let $E$ be an elliptic curve defined over $\QQ$. Then $E$ is modular, i.e., there exists a Hecke newform
 of weight two with Fourier expansion $f(q)= \sum_n b_nq^n$ such that for all primes $p$ of good
reduction
$$
a_p= 1 - N_p + p=b_p.
$$ 
\end{theorem}

Before commenting further on this theorem, we want to discuss an example. Let $E$ be the elliptic 
curve given by the equation
$$
y^2=x^3+1.
$$
Then we find the following values for $N_p$ and $a_p$:
$$
\begin{array}{r|rrrrrr rr }
p  & 2 & 3  & 5& 7 & 11 & 13 & 17 & 19\\
\hline
N_p & 3 & 4 & 6 & 12 & 12 & 12 & 18 & 12 \\
\hline
a_p & 0 & 0 & 0 & -4 & 0 & 2 & 0 & 8 
\end{array}
$$
In this case it is easy to describe the associated modular form explicitly. 
Recall that the Dedekind $\eta$-function is defined by 
$$
\eta(\tau)=\eta(q)= q^{\frac{1}{24}}\prod_{n \geq 1}(1-q^n).
$$
Then the form
$$
f(\tau)=\eta(6\tau)^4=q\prod_{n \geq 1}(1-q^{6n})^4= 1 - 4q^7 + 2 q^{13} + 8 q^{19} + \ldots
$$
is a weight two form of level $36$ and one can easily verify the condition $a_p=b_p$
for any primes within one's computing capacity. One can provide a complete proof for this fact using Theorem~\ref{theo_FSL} or the theory of CM-elliptic curves and CM-modular forms.

Some comments regarding Theorem~\ref{Thm_Wiles} are in order at this point. 
The Taniyama-Shimura-Weil conjecture also makes a prediction about the level
of the modular form. If the equation of $E$ is a \emph{global minimal Weierstrass equation} 
(this is a Weierstrass equation for $E$
whose discriminant $\Delta$ is divisible by a minimal power of $p$ for all primes $p$, 
see \cite[Chapter X.1]{Kn}), then the level of
$f$ is only divisible by the primes of bad reduction. 
More precisely, if $p \neq 2,3$, then the level is at most divisible by $p^2$, for $p=3$ it is at most 
divisible by $p^5$ and for $p=2$ by $p^8$ (see \cite{Ogg}). In fact, the level equals the conductor of $E$ (see the introduction of \cite{BCDT}.

Also, the statement that $E$ is modular, can be expressed in different, though equivalent, ways. One
of these is to say that an elliptic curve $E$ is modular, if there exists a non-constant morphism
$X_0(N) \to E$ defined over the rationals. This characterisation has no satisfactory known counterpart in the 
higher dimensional case. 

Yet another equivalent formulation, and one that can be generalised, 
is to say that the 
$L$-series of $E$ is that of a weight two form $f$. In Section~\ref{sub_Lfunctions} we have already 
defined an $L$-function for
varieties defined over the $\QQ$. In the case at hand we can make this explicit, including
a definition of the Euler factors associated to the bad primes:
we first assume that the elliptic curve $E$ is given by a global minimal
Weierstrass equation. Let $a_p=1-N_p+p$ for good primes and define
$$
a_p=\begin{cases}
1 \quad&\text{in the case of split multiplicative reduction} \\
-1 \quad&\text{in the case of non-split multiplicative reduction} \\
0 \quad&\text{in the case of additive reduction.} 
\end{cases}
$$ 
Then we set
$$
L(E,s)=\prod_{p\mid \Delta}\frac{1}{1-a_pp^{-s}}\prod_{p\nmid \Delta}\frac{1}{1-a_pp^{-s}+p^{1-2s}}.
$$
In Section~\ref{subsec_Galoisassmodform} we have also introduced the $L$-series of a normalised
Hecke eigenform.
Modularity of $E$ can then be rephrased as an 
equality of $L$-series
$$
L(E,s)=L(f,s)
$$
for a suitable Hecke eigenform of weight two.

\subsection{The theorem of Dieulefait and Manoharmayum}\label{sub_Dieulefait}

In view of Wiles' result it is natural to ask for generalisations to other classes of varieties.
Probably the most natural class to be considered for this question is that of Calabi-Yau varieties. 
A smooth projective variety is called a \emph{Calabi-Yau} variety if 
the canonical bundle is trivial, i.e.,
\begin{equation}\label{eqn_trivialcononical}
\omega_X = \cO_X
\end{equation}
and the following vanishing holds
\begin{equation}\label{eqn_vanishing}
h^i(X, \cO_X)=0 \, \mbox{ for } 0 < i <\dim X.
\end{equation}
For a curve, condition (\ref{eqn_vanishing}) is empty and we obtain the class of elliptic curves. Calabi-Yau 
varieties of dimension two are $\Kthree$ surfaces (see Section~\ref{sub_examplesK3} for a discussion 
of this case). Calabi-Yau threefolds have been investigated
very thoroughly. A lot of this interest comes from the fact that they play an essential
role in string theory. The easiest examples of Calabi-Yau threefolds are quintic 
hypersurfaces in $\PP^4$
and complete intersections of type $(2,4)$ or $(3,3)$ in $\PP^5$. For a discussion  
of Calabi-Yau varieties, which appear
as complete intersections in toric varieties, and their connection with mirror symmetry we 
refer the reader to \cite{Bat}, \cite{BaBo}.

If we look at Calabi-Yau varieties from the point of view of modularity, then it is natural to start
with examples whose middle cohomology is as simple as possible. By definition, we always have
$h^{3,0}(X)= h^{0,3}(X)=1$. A Calabi-Yau variety $X$ is called \emph{rigid} if $h^{2,1}(X)= h^{1,2}(X)=0$.
The name comes from the fact that this condition implies $H^1(X,T_X)=0$, and hence
that $X$ has no infinitesimal complex deformations. In this case the middle cohomology
$$
H^3(X)= H^{3,0}(X) \oplus H^{0,3}(X) \cong \CC^2
$$
is two-dimensional.

Now assume that $X$ is a rigid Calabi-Yau manifold, which is defined over $\QQ$. Then one can show
 that
the determinant of
Frobenius on $H^3_{\etal}(\overline{X},\QQ_l)$ is $p^3$. 
Hence
$$
P_{3,p}(t)= \det(1-t \Fr_p^*\mid H^3_{\etal}(\overline{X},\QQ_l))
= 1 - \trace((\Fr_p)^*\mid H^3_{\etal}(\overline{X},\QQ_l))t + p^3t^2.
$$
We shall again set
$$
a_p= \trace(\Fr_p^*\mid H^3_{\etal}(\overline{X},\QQ_l)).
$$
The $L$-series of (the middle cohomology of) $X$ is then of the form
$$
L(X,s)=(*) \prod_{p\in \mathcal{P}} 
\frac{1}{1-a_pp^{-s} + p^{3-2s}}
$$
where the factor $(*)$ denotes the Euler factors associated to the bad primes, and $\mathcal{P}$ is the product of good primes.
Comparing this to equation (\ref{eqn_Lseriesmodform}) leads one to expect a relationship
with a modular form of weight four. 
Indeed, modularity of rigid Calabi-Yau threefolds defined
over $\QQ$, had been conjectured for some time. This can be seen as a special form of the Fontaine-Mazur
conjecture \cite{FM} (see Subsection~\ref{sub_FM}) and has been stated explicitly by Saito and Yui~\cite{SY}. In recent years, many 
examples of geometrically interesting rigid Calabi-Yau varieties have been found and in all known cases 
it has been possible to determine the associated weight four cusp form (see also Section~\ref{sub_FSL}). 

A general modularity result has been proved by Dieulefait and Manoharmayum \cite{DM}, with further improvements
provided later by Dieulefait \cite{Di2}, \cite{Di3}.
\begin{theorem}\label{theo:DM}
Let $X$ be a rigid Calabi-Yau threefold defined over $\QQ$, and assume that one of the following
conditions holds:
\begin{enumerate}
\item $X$ has good reduction at $3$ and $7$ or
\item $X$ has good reduction at $5$ or
\item $X$ has good reduction at $3$ and the trace of $\Fr_3$ on 
$H^3_{\etal}(\overline{X},\QQ_l)$ is not
divisible by $3$.
\end{enumerate}
Then $X$ is modular. More precisely
$$
L(X,s) \circeq L(f,s)
$$
for some weight four modular form, where $\circeq$ means equality up to finitely many Euler factors. 
\end{theorem}

Very recently Dieulefait has informed us that the above theroem can be further improved. More precisely,
he can show that $X$ is modular if it has good reduction at $3$.   
The proof goes along the lines of \cite{Di3}.
This requires that a modularity lifting
result of Diamond, Flach and Guo can be extended to the case of weight $4$ and chracteristic $3$.
Indeed, Kisin \cite{Kis} has recently announced such an extension.

We remark that, as in the case of elliptic curves, one expects that the level of the modular form is 
only divisible by the primes of
bad reduction. 
For the question of determining the level of $f$ see also \cite{Di1} and Section~\ref{Sec_open}.
Very little is known about the Euler factors associated to the primes of bad reduction.

The above result is purely an existence result and does not provide a method to determine the form
$f$ explicitly. There is, however, a method due to Faltings, Serre and Livn\'e, which is very
effective, if one wants to prove that a candidate modular form is indeed the right one. 
We shall discuss this
method in the next section.

Finally, we would like to remark that a number of non-rigid Calabi-Yau varieties have been found,
for which modularity has been established. Most of the examples known up to date are of two types. They
are either of Kummer type or they contain (many) elliptic ruled surfaces (or both).
In both cases the geometry of the variety is such that the middle cohomology breaks up into
two-dimensional pieces.
The first type
of example was found by Livn\'e and Yui \cite{LY}, for Calabi-Yau varieties containing elliptic ruled 
surfaces see \cite{HV1}  and \cite{HV2}. We shall discuss examples, mostly rigid, in more detail in Section~\ref{sec_examples}.

\subsection{The method of Faltings-Serre-Livn\'e}\label{sub_FSL}

Given a rigid Calabi-Yau threefold, defined over $\QQ$, one often wants to determine the 
associated modular form explicitly. Generally, it is not hard to find a suitable candidate for
such a form $f$: counting points, one can determine the first few Fourier coefficients of $f$ and
then compare these to existing lists of Hecke eigenforms such as \cite{St}. In order
to actually prove that the $L$-series of the variety and that of the modular form coincide,
one can try to establish the equality of (the semi-sim\-pli\-fi\-ca\-tions of) the corresponding two-dimensional
Galois representations (over $\QQ_2$).
In order to prove that two Galois representations $\rho_1,\rho_2$ have isomorphic semi-sim\-pli\-fi\-ca\-tions, it suffices to prove that for all primes $p$ in a certain finite set of primes $T$ the traces of $\rho_1(\Fr_p)$ and $\rho_2(\Fr_p)$ are equal. This result is due to Faltings \cite{Fa}. In the sequel one of the representations is the Galois representation on the cohomology, the other representation is of the from $\rho_f$, with $f$ a modular form. It is known that $\rho_f$ is already simple, see \cite[Theorem 2.3]{Ri}.

If the representations have even traces
then one can determine the set $T$ effectively. This  method is  due to Serre and was recast by Livn\'e \cite{L1}. 
Otherwise one can use an approach of Serre \cite{Se}, which, however can only be made explicit in special
cases.

\begin{theorem}\label{theo_FSL}
Let $\rho_1, \rho_2$ be two continuous two-dimensional two-adic
representations of $\Gal(\overline{\QQ}/\QQ)$, unramified outside a
finite set $S$ of prime numbers. Let $\QQ_S$ be the compositium of all
quadratic extensions of $\QQ$, which are unramified outside $S$ and let $T$ be a
set of primes, disjoint from $S$, such that
$\Gal(\QQ_S/\QQ)=\{\Fr_p\mid_{\QQ_S}; p\in T\}$.
Suppose that
\begin{enumerate}
\item 
$\trace\rho_1(\Fr_p)=\trace\rho_2(\Fr_p)$ for all $p\in T$,
\item 
$\det \rho_1(\Fr_p)=\det \rho_2(\Fr_p)$ for all $p\in T$,
\item
$\trace \rho_1\equiv \trace \rho_2 \equiv 0 \bmod 2$
 and $\det \rho_1\equiv \det \rho_2 \bmod  2$. 
\end{enumerate}
Then $\rho_1$ and $\rho_2$ have isomorphic semisimplifications, and hence
$L(\rho_1, s)\circeq L(\rho_2, s)$. In particular, the good Euler factors of
$\rho_1$ and $\rho_2$ coincide.
\end{theorem}

The gist of this theorem is, that, provided the parity of trace and determinant coincide, it is
sufficient to check only the equality of the traces and determinants for a \emph{finite} number of explicitly known primes, 
in order to establish that the semi-sim\-pli\-fi\-ca\-tions are isomorphic. 

To check the parity of trace and 
determinant is in practice often easy. The main point is to determine a suitable set $T$ of primes.
Given this, the calculation of the traces $a_p$ for $p \in T$ is often straightforward. Of course, the
set $T$ depends on the bad primes of $X$. But in most concretely given examples, the set $T$ is quite small.

If the traces are not even, then one can apply Serre's approach, provided that the 
mod $2$-reductions $\bar\rho_1$ and $\bar\rho_2$ are absolutely irreducible and isomorphic. 
This means that their kernels cut out the same $S_3$-extension $K$ of $\QQ$. 
In this situation, one 
has to compute all  extensions $L$ of $K$ such that $\Gal(L/\QQ)$ is either  $S_4$ or $S_3\times C_2$ and which are unramified outside 
the ramification loci of the representations $\rho_i$. 
Let $T$ be a set of primes such that for any such extension $L$ there is a prime $p\in T$ with $\Fr_p$ of 
maximal order in Gal$(L/\QQ)$. Then, proving the isomorphism of the Galois representations amounts to 
checking equality of the traces at the primes of $T$.

In practice, this method can only be used if there are only a few small bad primes, since the available 
tables of such number fields are limited (cf. \cite{Jo}). 
We will sketch this approach for the Schoen quintic in Section~\ref{subsub_Schoenquintic}. 
In the context of non-rigid Calabi threefolds, applications of this method can also be found in \cite{Sch2}.

\section{Examples}\label{sec_examples}

In this chapter we shall discuss various examples which have been found in recent years. For an
extensive survey of modular Calabi-Yau threefolds we refer the reader to the book of Meyer \cite{Me}.

\subsection{$\Kthree$ surfaces}\label{sub_examplesK3}

In Sections~\ref{sub-Wiles} and~\ref{sub_Dieulefait} we discussed elliptic curves and (rigid) 
Calabi-Yau threefolds. 
The reader may ask why we did not consider $\Kthree$ surfaces before proceeding to dimension three.
The reason is that in even dimension the middle cohomology will always contain algebraic classes, which
changes the situation. If $S$ is any $\Kthree$ surface, then the middle cohomology has dimension $22$, more
precisely
$$
H^2(S,\ZZ) = 3U + 2E_8(-1)
$$
where $U$ denotes the hyperbolic plane, i.e. the free lattice of rank $2$, equipped with the intersection 
form given by the matrix $\begin{pmatrix} 0 & 1\\ 1&0 \end{pmatrix}$ and $E_8$ is the unique 
positive definite, even, unimodular 
lattice of rank $8$. The notation $E_8(-1)$ indicates that we take the negative 
form. The N\'eron-Severi group of $S$ is the group of divisors modulo numerical equivalence or, equivalently,
$$
\NS(S)= H^{1,1}(S) \cap H^2(S,\ZZ).
$$
Its rank $\rho(S)$ is called the \emph{Picard number} of $S$. This can vary (in characteristic $0$)
from $0$, in which case the $\Kthree$ surface is not algebraic, to $20$. $\Kthree$ surfaces $S$ with maximal
Picard number $\rho(S)=20$ are called \emph{singular} $\Kthree$ surfaces (a common but misleading
terminology which does not mean that the surface has singularities). Singular $\Kthree$ surfaces have no moduli,
but they form an everywhere dense set (in the analytic topology) in the period domain of 
$\Kthree$ surfaces. The orthogonal complement of  
$\NS(S)$ in $H^2(S,\ZZ)$ is denoted by $T(S)$ and is called the \emph{transcendental lattice}
of the surface $S$. If $S$ is a singular $\Kthree$ surface, this is a positive definite even 
lattice of rank $2$. Shioda and Inose have proved that the map which associates to a 
singular $\Kthree$ surface $S$
the oriented lattice $T(S)$ (where the orientation is given by the $2$-form on $S$) 
defines a bijection between the sets of singular $\Kthree$ surfaces modulo isomorphism
on the one hand
and positive definite, integral binary quadratic forms modulo $\SL(2,\ZZ)$ on the other hand \cite{SI}.
One also refers to the discriminant of the lattice $T(S)$ as the \emph{discriminant} of the surface $S$. 
The injectivity of the above map is a consequence of the Torelli theorem. To prove surjectivity Shioda and
Inose  proceed as follows. Fix a positive definite even lattice of rank $2$ and discriminant $d>0$.
Recall that such a lattice
determines uniquely a \emph{reduced} quadratic form
$$
Q= \begin{pmatrix} 
2a & b\\ b& 2c 
\end{pmatrix}
$$
where $\det(Q)=d$ and $a>0, -a <b\leq a \leq c$ (and $b \geq 0$ if
$a=c$). Set
$$
\tau=\frac{-b + \sqrt{-d}}{2a} \mbox{ and } \tau'=\frac{b + \sqrt{-d}}{2}.
$$
By construction, the elliptic curves $E_{\tau}$ and $E_{\tau'}$ have complex multiplication
and are isogenous. The abelian surface $A=E_{\tau} \times E_{\tau'}$ has Picard number $\rho(A)=4$
and the intersection form on its transcendental lattice $T(A)$ is given by $Q$, see \cite{MiSh}. For the Kummer surface
$\Km(A)$ one has $T(\Km(A))=T(A)(2)$, i.e. the form defined by $Q$, but multiplied by $2$. 
In order to find a $\Kthree$ surface with the form given by $Q$, Shioda and Inose exhibit an elliptic fibration 
on $\Km(A)$ and construct a suitable double cover, resulting in another $\Kthree$ surface $S$ for which
$T(S)=T(A)$.
The surfaces $S$ and $\Km(A)$ are related by a \emph{Nikulin involution}. More precisely, 
$S$ admits an involution $\iota$ with exactly $8$ fixed points, leaving the non-degenerate
$2$-form on $S$ invariant, such that $\Km(A)$ is the minimal resolution of the quotient $S/\iota$.

We shall now turn to the question of modularity of $\Kthree$ surfaces. Note that we are now mostly considering
varieties defined over a number field $K$ rather than over $\QQ$. We can, however, still define an $L$-series
for a variety $X$ defined over $K$. In order to do this, we define for every prime ideal $\mathfrak{p} \in \Spec \cO_K$ and 
each integer $i$ the polynomial 
$$
P_{i,\mathfrak{p}}(t)= \det(1- t\Fr_{\mathfrak{p}}^*\mid H^i_{\etal}(\overline{X}_{\mathfrak{p}},\QQ_l)).
$$ 
If the dimension of $X$ is $d$, we set
$$
L(X/K,s) = (*)\prod_{\mathfrak{p}\in \mathcal{P}} \frac{1}{P_{d,\mathfrak{p}}((\Nm\mathfrak{p})^{-s})},
$$
where $\Nm$ denotes the norm endomorphism from $K$ to $\QQ$ and $\mathcal{P}$ is the set of primes of good reduction.
The symbol $(*)$ stands for suitably defined
Euler factors at the bad primes. In the case $K=\QQ$ this coincides with our previous definition.

Since very little is known about the modularity of $\Kthree$ surfaces which are not singular,
we shall restrict ourselves to the singular case.
Recall that a modular form $f$ with Fourier expansion
$f = \sum_n a_n q^n$ has complex multiplication
by a Dirichlet character $\chi$ if $f=f \otimes \chi$ where
$$
f \otimes \chi = \sum_n \chi(n) a_n q^n.
$$
The character $\chi$ is necessarily quadratic. We say that $f$ has CM by the imaginary quadratic field 
$K$, if $f$ has CM by the corresponding quadratic character. CM forms are closely related to 
Hecke Gr\"ossenscharacters: by a result of Ribet, every newform with CM 
comes from a Gr\"ossencharacter of an imaginary quadratic field $K$. For a treatment of the connection between
CM forms and Hecke Gr\"ossencharacters see~\cite{Ri}.

Let $S$ be a singular $\Kthree$ surface of discriminant $d$ and $K_0= \QQ(\sqrt{-d})$. 
One can find a finite extension $K$ of $K_0$, such that $S$ has a model over $K$
and that, moreover, the N\'eron-Severi group of $S$ is generated by divisors which are defined 
over $K$. Shioda and Inose \cite[Theorem 6]{SI} proved 
\begin{theorem}
Let $S$ and $K$ as above. Then
$$
L(S/K,s) \circeq \zeta_K(s-1)^{20} L(s,\psi^2)L(s,\overline{\psi}^2)
$$ 
where $\zeta_K(s)$ is the Dedekind zeta function of the field $K$ and
$\psi$ is the Gr\"ossencharacter of a model of the elliptic curve $E_{\tau}$ over $K$ which is associated to $S$.
Here $\circeq$ denotes equality up to finitely many factors.
\end{theorem}
 
We shall now turn to the more special case that the singular $\Kthree$ surface $S$ is defined over $\QQ$.
Then the transcendental lattice $T(S)$ is a two-dimensional $\Gal(\overline{\QQ}/\QQ)$-module, which defines
an $L$-series $L(T(S),s)$. This case was treated by Livn\'e \cite[Theorem 1.3]{L2}.
\begin{theorem}  
Let $S$ be a singular $\Kthree$ surface of discriminant $d$, defined over $\QQ$. Then there exists 
a weight three modular form $f_3$ with complex multiplication by $\QQ(\sqrt{-d})$ such that 
$$
L(T(S),s) \circeq L(f_3,s).
$$
\end{theorem}

If, moreover, the N\'eron-Severi group $\NS(S)$ is generated by divisors defined over $\QQ$,
then one has in fact that 
$$
L(S,s) \circeq \zeta(s-1)^{20} L(f_3,s).
$$

The level of the form $f_3$ equals the discriminant of $S$ up to squares composed of the bad primes \cite[Proposition 13.1]{Sch3}.

An example where the latter situation occurs is the universal
elliptic curve $\mathcal{E}_1(7)$ with a point of order $7$. An affine model of this surface is given by
$$
 y^2+(1+t-t^2)xy+(t^2-t^3)y=x^3+(t^2-t^3)x^2.
$$
In this case, $\NS(\mathcal{E}_1(7))$ is generated by divisors defined over $\QQ$, and 
the associated modular form is $f_3= (\eta(\tau)\eta(7\tau))^3$, where $\eta$ is the 
Dedekind $\eta$-function.
For details of this example see \cite{HV3}.
Another example of this kind was considered by Sch\"utt and Top \cite{ST}, who investigated the 
extremal elliptic fibration $\pi: S \to \PP^1$ with singular fibres $I_{19},I_1,I_1,I_1,I_1,I_1$.

At this point we would like to mention the following result, due to Shafarevich \cite{Sha}. 
\begin{theorem}\label{theo_Shafarevich}
Given an integer $n$, there is a finite set $\cS_n$ of $\SL_2(\ZZ)$-equi\-va\-len\-ce classes of even, positive definite, binary quadratic
forms, such that, if $S$ is a singular $\Kthree$ surface admitting a model over a number field of degree 
at most $n$ then $T(S) \in \cS_n$.
\end{theorem}
In particular, this result says that there exist only finitely many $\CC$-isomorphism classes of 
singular $\Kthree$ surfaces which can be defined over $\QQ$.

It is a natural question to ask which weight three newforms with rational coefficients can be realized by singular
$\Kthree$ surfaces defined over $\QQ$. In view of Shafarevich's result, this number is finite (up to twisting).
Motivated by this observation Sch\"utt \cite{Sch3} proved the following finiteness result for CM forms:
\begin{theorem}\label{theo_CMforms}
Assume the generalised Riemann hypothesis (GRH). Then, for given weight, there are (up to twisting) only finitely 
many CM newforms with rational coefficients. For weights $2,3$ and $4$ this holds unconditionally.
\end{theorem}

For weight two this is, of course, classical theory of CM elliptic curves. The fact that this theorem
is true without the assumption of GRH for weights three and four follows from Weinberger's result \cite{We} of 
the finiteness
of imaginary quadratic fields with exponent $2$ or $3$. Sch\"utt has determined lists of such 
forms for weights three and four. For weight three this list is complete under the assumption of GRH and
without this assumption the list is complete up to at most one form of level at least $2 \cdot 10^{11}$. For weight four the list is complete up to
level $10^{10}$.
So far, he has found singular $\Kthree$ surfaces over $\QQ$ for $24$  out of the $65$ weight three forms.
\subsection{Rigid Calabi Yau threefolds}\label{sub_rigidCY}

In this section we will introduce some of the rigid Calabi-Yau threefolds over $\QQ$ for which modularity 
has been proven. For a complete treatment, the reader is referred to the book of Meyer \cite{Me}. 
Nevertheless, we tried to include the most important and instructive examples. 
Note that for all of them modularity statements can also be derived from Theorem~\ref{theo:DM}
(although this theorem does not give the explicit modular form).

\subsubsection{The Schoen quintic}\label{subsub_Schoenquintic}

To our knowledge the first rigid Calabi-Yau threefold over $\QQ$ for which modularity was
proved is the Schoen quintic $X$ in $\PP^4$:
\[
X: x_0^5+x_1^5+x_2^5+x_3^5+x_4^5-x_0x_1x_2x_3x_4=0.
\]
Note that $X$ appears as a (singular) member of the one-parameter family $\cM$ in 
Section~\ref{sub_zetamirror}. Schoen \cite{Schoen1} showed that $X$ is rigid and derived the 
modularity of $X$ explicitly:

\begin{theorem}
$X$ is modular and its $L$-series coincides with the Mellin transform of a 
normalised newform of weight four and level $25$.
\end{theorem}

Since the newform does not have even coefficients, the proof of the theorem requires 
Serre's approach (cf. Section~\ref{sub_FSL}). In this case, this is easy, since 
by \cite{Jo} (cf. \cite[Proposition 4.10]{L1}) there is a 
unique Galois extension $K/\QQ$ with Galois group $S_3$ or $C_3$ which is unramified outside $2$ and $5$.
By \cite{Jo}, there are only six extensions with Galois group $S_4$ or $S_3\times C_2$  
which have to be considered. 
One then finds that it suffices to check the primes in the set $T=\{3,7,11,13\}$.

\subsubsection{Fibre products of rational elliptic surfaces}\label{subsub_fibre}

A very useful method to construct Calabi-Yau threefolds goes back to Schoen. 
In \cite{Schoen2} he considers fibre products of two regular semi-stable rational elliptic surfaces with a section. 
Due to the semi-stability, such a fibre product has only ordinary double points as singularities, and hence 
a small resolution exists in the analytic category. Such a resolution has all properties of a Calabi-Yau
variety, but may not be projective. 

Sometimes, it is possible to construct a small projective resolution by successively blowing up suitable 
subvarieties. 
In general, there are two major classes of fibre products where enough subvarieties are available: 
\begin{proposition}[Schoen {\cite[Lemma 3.1]{Schoen2}}]\label{prop_Schoen} For $i=1,2$ let $\pi_i: Y_i\rightarrow \PP^1$ be a semi-stable rational elliptic surface with a section.

Let $X$ be the fibre product $Y_1\times_{\PP^1} Y_2$.
Then $X$ has a projective small resolution $\hat{X}$ if
\begin{enumerate}
\item all singularities of $X$ lie in fibres of type $I_m\times I_n$ with $m, n>1$, or
\item $Y_1=Y_2$ and $\pi_1=\pi_2$ .
\end{enumerate}
\end{proposition}

Schoen also investigated when a resolution of a fibre product is rigid. 
This is independent of its projectivity. He describes four cases where this happens. 
In this section, we will see examples for two of them.

Let $X$ be the self-fibre product of a semi-stable rational elliptic surface with the minimal number of 
four singular fibres. Then $X$ has a projective small resolution $\hat X$ by the proposition above. This is 
rigid due to \cite[Proposition 7.1, (i)]{Schoen2}. There are six such surfaces up to isomorphism over $\overline{\QQ}$, as determined by 
Beauville \cite{B}. 
All of them are modular in the sense of Shioda \cite{ShiMod}, i.e., they are elliptic surfaces associated with a subgroup of $\SL(2,\ZZ)$ of finite index. We give the corresponding congruence subgroups in  Table~\ref{Tab_Beau}, where we also include the types of the singular fibres.

The fact that the elliptic surfaces can be defined over $\QQ$,  and that $X$ is a self-fibre product imply that the small resolution $\hat X$ of $X$ is defined over $\QQ$. 
Hence, we expect it to be modular. 

There are three alternative proofs of modularity: one uses the Eichler-Shimura 
isomorphism and the fact that the respective spaces  of newforms are 
one-di\-men\-sio\-nal.
This has been worked out in \cite{SY}. The next approach consists of counting 
points and 
applying the Faltings-Serre-Livn\'e method as sketched in Section~\ref{sub_FSL}. 
This was done by Verrill in \cite[Appendix]{Y}. We note that due to the 
self-fibre product structure, 
the number of points on the reductions and, consequently, the traces are 
always even, so Theorem~\ref{theo_FSL} 
applies. Then one can use Theorem~\ref{theo:DM} and then try to proceed as 
in \cite{Di1}.

Finally, take for every $\overline{\QQ}$-isomorphism class the elliptic surface given by the equations listed in \cite{B}.  Table~\ref{Tab_Beau} gives the weight four newforms associated to the self-fibre 
products and their respective levels. 

We also remark that these newforms are the only ones of weight four which 
can be written as $\eta$-products (see \cite[Cor. 2]{Ma}).

\begin{table}
$$
\begin{array}{|c|c|c|c|}
\hline
\Gamma & \text{singular fibres} & \text{weight $4$ form} & \text{level}\\
\hline
\Gamma(3) & I_3, I_3, I_3, I_3 & \eta(3\tau)^8 & 9\\
\hline
\Gamma_1(4)\cap\Gamma(2) & I_4, I_4, I_2, I_2 & \eta(2\tau)^4\eta(4\tau)^4 & 8\\
\hline 
\Gamma_1(5) & I_5, I_5, I_1, I_1 & \eta(\tau)^4\eta(5\tau)^4 & 5\\
\hline
\Gamma_1(6) & I_6, I_3, I_2, I_1 & \eta(\tau)^2\eta(2\tau)^2\eta(3\tau)^2\eta(6\tau)^2 & 6\\
\hline
\Gamma_0(8)\cap\Gamma_1(4) & I_8, I_2, I_1, I_1 & \eta(4\tau)^{16}\eta(2\tau)^{-4}\eta(8\tau)^{-4} & 16\\
\hline
\Gamma_0(9)\cap\Gamma_1(3) & I_9, I_1, I_1, I_1 & \eta(3\tau)^8 & 9\\
\hline
\end{array}
$$
\caption{Beauville surfaces}\label{Tab_Beau}
\end{table}

We shall give one further class of examples which was studied in \cite{Sch1}. 
As above, we shall work with the universal elliptic curve with a point of order $6$, which we will
denote by $\mathcal{E}_1(6)$. (This is the elliptic surface associated to $\Gamma_1(6)$.) 
We take $\pi_2=\sigma\pi_1$ with $\sigma$ an automorphism of $\PP^1$ permuting the three cusps with reducible fibres and fixing the remaining cusp.

This leads to case 1 of Proposition~\ref{prop_Schoen}.
The projective small resolution $\hat X$ is rigid by \cite[Proposition 7.1, (ii)]{Schoen2}. There are five 
such automorphisms of $\PP^1$ which are all defined over $\QQ$. 
The resulting rigid Calabi-Yau threefolds 
can be proved to be modular and one determines the corresponding newforms using Theorem~\ref{theo_FSL}. The associated newforms 
are the forms of  level $10, 17, 21,$ and twice the same form of level $73$.

\subsubsection{Further examples}\label{subsub_further}
We shall briefly comment on other modular rigid Calabi-Yau threefolds. 
To our knowledge, the first examples which were published  after the Schoen quintic were 
obtained by Werner and Van Geemen \cite{WvG}, based on work of Hirzebruch. 
Their paper pioneered an important method to determine the third Betti number $b_3(X)$ for a 
Calabi-Yau threefold $X$ over $\QQ$ arithmetically. This method
requires the knowledge of the Euler number 
$e(X)$ and of a good prime $p$ such that $\Fr_p$ operates trivially
on the divisors of $X$. Then the Riemann hypothesis gives a bound in terms of $b_3(X)$  
on the trace of $\Fr_p$ on $H^3_{\etal}(X,\QQ_l)$. This bound
 sometimes enables one to deduce the value of $b_3(X)$ after 
counting points. Using this approach and Theorem~\ref{theo_FSL}, Werner and van Geemen proved modularity 
for four rigid Calabi-Yau threefolds over $\QQ$. They also determined the associated newforms, 
finding the forms of level $8$ and $9$ from above and two further forms of level $12$ and $50$ respectively. 

Further, we would  like to mention two constructions which produced several modular 
Calabi-Yau threefolds. The first of them uses toroidal geometry involving root lattices. 
For the root lattices $A_1^3$ and $A_3$, this was worked out by Verrill in \cite{Ve}. The corresponding newforms of 
level $6$ and $8$ appear in Table~\ref{Tab_Beau}. For $A_4$, a family of Calabi-Yau threefolds was investigated by 
Hulek and Verrill in \cite{HV1}. This family contains a number of (rigid and non-rigid) 
members for which they prove modularity.

Another approach uses double octics, as studied by Cynk and Meyer 
in \cite{CM}, \cite{CM3}. In this case one starts with an octic surface $S$ in $\PP^3$ which has mild 
singularities. Many such double covers of $\PP^3$ 
admit a desingularisation $X$, which is Calabi-Yau. In practice one uses for $S$ highly reducible surfaces,
such as a union of $8$ planes (such configurations are also 
called \emph{octic arrangements}).  Such an $X$ is called a double octic. 
Often it is possible to describe the deformations of a double octic $X$ explicitly in terms of the 
surface $S$. Using this, Cynk and Meyer found further rigid Calabi-Yau threefolds over $\QQ$ and the 
corresponding newforms in \cite{CM}.

One advantage of the constructions involving fibre products, root lattices or double octics is that 
they also produce non-rigid modular Calabi-Yau threefolds. Here, modularity is understood as the 
splitting of $H^3_{\etal}(X,\QQ_l)$ into two-dimensional Galois representations, all of 
which can be interpreted in terms of modular forms. In practice, this is often achieved by identifying 
enough elliptic ruled surfaces in $X$ which contribute to $H^3_{\etal}(X,\QQ_l)$. 
This method is explained in \cite{HV2}. For applications, the reader is also referred to 
 \cite{HV1}, \cite{Sch2} and \cite{CM3}.

Apart from this strict notion of modularity, the Langlands' programme hardly seems to be within reach in 
more generality at the moment. In \cite{CS}, Consani-Scholten investigate a Calabi-Yau threefold $X$ 
over $\QQ$ with $b_3(X)=4$. They prove that the Galois representation on 
$H^3_{\etal}(X,\QQ_l)$ is induced from a two-dimensional representation of 
$\Gal(\overline\QQ/\QQ(\sqrt{5}))$ and conjecture that this should come from a Hilbert modular form over 
$\QQ(\sqrt{5})$ of weight $(2,4)$ and level $5$. Recently Yi \cite{Yi} has announced a proof of this
conjecture.

\subsection{Higher-dimensional examples}\label{sub_higherdim}

A number of examples of higher-di\-men\-sio\-nal modular Calabi-Yau manifolds were recently constructed in
\cite{CH}. Most of these examples arise as quotients of products of lower dimensional Calabi-Yau 
manifolds by finite groups. The starting point is the Kummer construction, which was first used 
by Borcea \cite{Bor} and Voisin \cite{Vo}, who constructed Calabi-Yau threefolds as (birational) quotients
of a product of a $\Kthree$ surface admitting an elliptic fibration and an elliptic curve, by the diagonal
involution. This construction works in all dimensions: let $Y$ be a projective manifold with 
$H^q(Y,\cO_Y)=0$ for $q >0$ and $D \in \mid  -2K_Y\mid $ be a smooth divisor. Then $D$ defines a double cover
$\pi: X \to Y$, branched along $D$ with $X$ a Calabi-Yau manifold. Now assume that we have two 
coverings $\pi_i: X_i \to  Y_i$ as above. Then the group $\ZZ/2\ZZ \times \ZZ/2\ZZ$ acts on the product
$X_1 \times X_2$. 

\begin{proposition}\label{prop_Kummer}
Under the above assumptions, the quotient of the product $X_{1}\times X_{2}$ by the diagonal
involution admits a crepant resolution $X$, 
which is Calabi-Yau. Moreover, there is a double cover $X \to Y$, branched
along a smooth divisor $D$, with $H^{q}(Y,\cO_{Y})=0$ for $q>0$.  
\end{proposition}

The last remark of this proposition allows one to use the covering $\pi: X \to Y$ to repeat this 
process inductively. Starting with  low-dimensional modular Calabi-Yau varieties, one can in this way easily construct 
many higher dimensional Calabi-Yau varieties. The middle cohomology of these varieties grows with 
the dimension, but splits in many cases into two-dimensional pieces, which allows to prove modularity. 
Here we give one example. As the first factor we choose the rigid
Calabi-Yau threefold $X_{3}$, constructed as a resolution of
singularities of the double covering of $\PP^3$, branched along the following
arrangement of eight planes
$$
xt(x-z-t)(x-z+t)y(y+z-t)(y+z+t)(y+2z)=0.
$$
For details of this example see \cite[Octic Arr. No. $19$]{Me}. 
As the second factor we take the $\Kthree$ surface $S$ which is
obtained as a desingularisation of the double sextic 
branched along the following arrangement of six
lines
$$
xy(x+y+z)(x+y-z)(x-y+z)(x-y-z)=0.
$$
Performing the Kummer construction just described, we obtain
a smooth Calabi-Yau fivefold
$X_{5}$. The Hodge groups of $X_{5}$ are the invariant part of
the Hodge groups of $X_{3}\times S$. 
{From} this it follows easily that $b_{1}(X_{5})=b_{3}(X_{5})=0$ and $b_{5}(X_4)=4$. More precisely
$$
H^{5}(X_{5})\cong H^{3}(X_{3})\otimes T(S)
$$
where $T(S)$ is the transcendental lattice of $S$.
It is known that 
$$
L(T(S),s) \circeq L(g_3,s), \quad L(X_3,s) \circeq L(g_4,s)
$$
where $g_3$ and $g_4$ are the unique weight three and weight four Hecke eigenforms of level $16$ and $32$ 
respectively, with complex multiplication. 

Both of these forms can be derived from the unique weight two level $32$ newform
$$
g_2(q)=\eta(q^{8})^2\eta(q^{4})^2
$$
by taking the second, respectively,  third power of the corresponding Gr\"os\-sen\-cha\-rac\-ter $\psi$ of conductor $2+2i$. 
The $\Gal(\overline{\QQ}/\QQ)$-module $H^5_{\etal}(\overline{X},\QQ_l)$
is the tensor product of the $\Gal(\overline{\QQ}/\QQ)$-module $H^3_{\etal}(\overline{X},\QQ_l)$ and of the 
transcendental lattice $T(S)$ of $S$. Its semi-sim\-pli\-fi\-ca\-tion splits into two two-dimensional pieces and
one finds for the $L$-series 
$$
L(X_{5},s) \circeq L(g_{4}\otimes g_{3},s) = L(g_6,s)L(g_2,s-2)
$$
where $g_2$ is as above, and $g_{6}$ is a level $32$ cusp form of weight $6$, 
which can be derived from $g_2$ by taking the fifth power of the Gr\"os\-sen\-cha\-rac\-ter $\psi$.

If one wants to find examples of higher dimensional Calabi-Yau manifolds which are the analogue of
rigid Calabi-Yau varieties, i.e., have two-dimensional middle cohomology, then one has to 
divide by more complicated groups. One example of such a construction is the following. Consider
the elliptic curve 
$$
E: y^2=x^3-D
$$
for some positive integer $D$.
Then $x \mapsto \rho x$, with $\rho$ a third root of unity, defines an automorphism $\eta$ of order $3$ on $E$.
\begin{theorem}\label{theo_higherdim_E} 
Let $n>1$ be an integer and let $X'$ be the quotient of the cartesian product $E^{n}$ by 
the group 
$$
\{(\eta^{a_{1}}\times\dots\times \eta^{a_{n}})\in
\operatorname{Aut}(E^{n}); a_{1}+\dots+a_{n}\equiv0\bmod3\}.
$$
Then $X'$ has a smooth model $X$ which is a Calabi-Yau
manifold and the dimension of $H^n_{\etal}(\overline{X},\QQ_\ell)$ equals  two if $n$ is odd, respectively the dimenson of $T(X)\otimes \QQ_l$ equals two if $n$ is even,
where $T(X)$ is the transcendental part of the middle cohomology.

Moreover, $X$ is defined over $\QQ$ and $L(H^n_{\etal}(\overline{X},\QQ_l),s)\circeq L(g_{n+1},s)$, respectively
$L(T(X)\otimes \QQ_l,s)\circeq L(g_{n+1},s)$, where
$g_{n+1}$ is the weight $n+1$ cusp form with complex multiplication in
$\QQ(\sqrt{-3})$, associated to the $n$-th power of the Gr\"ossencharacter of $E$. 
\end{theorem}

One can also attempt to generalise the construction of double covers of $\PP^n$, branched along a configuration of hyperplanes,
to higher dimension. In general, however, this becomes very difficult. The combinatorial problems in order to control
the singularities of, say, a double cover of $\PP^5$ branched along  $12$ hyperplanes become very hard. In \cite{CH} one example was studied, namely the
double cover $X$ of $\PP^5$ given by
$$
w^2= x(x-u)(x-v)y(y-u)(y-v)z(z-u)(z-v)t(t-u)(t-v).
$$ 
This example was first studied by Ahlgren, who counted the number of points modulo $p$ on the 
affine part of $X$, relating these numbers to a modular form. The following was proved
in \cite[Theorem 5.12]{CH}.

\begin{theorem} \label{theo_Ahlgren}
The variety $X$ has a smooth model $Z$, defined over $\QQ$, 
which is a Calabi-Yau fivefold with Betti numbers 
$b_{1}(Z)=b_{3}(Z)=0$, $b_{5}(Z)=2$. More precisely, $h^{5,0}=h^{0,5}=1$.
One has 
$$
L(Z,s) \circeq L(f_6,s)
$$
where $f_6(q)=\eta^{12}(q^2)$ is the unique normalised cusp form of level four and weight six. 
\end{theorem}   

\section{Zeta Functions and Mirror Symmetry}\label{sec_Mirror}
\subsection{ $p$-adic cohomology}\label{sub_padiccohomology}
In practice, it turns out to be hard to work with \'etale cohomology in an explicit way (at least 
in dimension two and higher). In this subsection we discuss another approach to calculate the 
zeta-function of a smooth variety over a finite field $\FF_q$. This gives enough information to 
obtain the characteristic polynomial of $\Fr_p$ on $H^i_{\etal}(X,\QQ_{l})$ for a prime of good reduction. 
We use cohomology with $p$-adic coefficients instead of $l$-adic coefficients. 
Cohomology theories with $p$-adic coefficients occured, for example, in Dwork's proof \cite{Dw} of the 
rationality of the zeta-function.

There exist several $p$-adic cohomology theories which are rich enough to prove the Weil conjectures, 
and therefore are called `Weil cohomologies'. These theories have a Lefschetz trace formula. 

We discuss an easy variant, called Mon\-sky-Wash\-nit\-zer cohomology, which strictly speaking is not a Weil cohomology. For a survey on this theory we 
refer to \cite{vdP}.

Let $p$ be a prime number, $r$ a positive integer. Set $q=p^r$.
Suppose $V/\FF_q$ is a hypersurface in $\PP^n$ and let $X$ be the complement of $V$ in $\PP^n$. 
Then $X$ is an affine variety, say $X=\Spec R'$. Consider the ring $\ZZ_q$,  the ring of 
integers of $\QQ_q$, which is the unique unramified extension of $\QQ_p$ of 
degree $r$. Let $\pi$ be the maximal ideal of $\ZZ_q$.
The ring $R'$ can be `lifted' to $\ZZ_q$, i.e., there is an 
algebra $\hat{R}:=\ZZ_q[X_0,\dots,X_m]/(f_1,\dots,f_k)$ such that $\hat{R}/\pi \hat{R}\cong R'$. 

Let $R$ be the `overconvergent completion' of $\hat R$, that is the ring 
\[\frac{\{ g \in \ZZ_q[[X_0,\dots,X_m]] ; g \mbox{ converges on an open disc of radius  more than } 1\}}
{(f_1,\dots, f_k)}.\]

This definition enables us to have a good theory of differential forms on $X$: consider the complex 
of differential algebras
\[ 0 \rightarrow \Omega^0_R\rightarrow \Omega^1_R\rightarrow \Omega^2_R\rightarrow \dots\rightarrow 
\Omega^n_R\rightarrow 0.\]
Then we define $H^i_{MW}(X,\QQ_q)$ to be the $i$-th cohomology group of this complex tensored 
with $\QQ_q$. It turns out that these groups are finite-dimensional.
The reason to work with (overconvergent) $p$-adic algebras is that they admit lifts of Frobenius, i.e., there exists a morphism $\Fr_q: R\rightarrow R$ such that  the reduction $\Fr_q:R/\pi R\rightarrow R/\pi R$ is the geometric Frobenius map on $X$.

One can show that 
\[ \sum_i (-1)^i \trace( q^n\Fr_q^{-1} \mid H_{MW}^i(X,\QQ_q))=\# X(\FF_q).\]
The appearance of $q^n\Fr_q^{-1}$ is comes from the fact that Mon\-sky-Wash\-nit\-zer cohomology 
is the Poincar\'e dual of rigid cohomology with compact support. The latter theory has a more usual 
Lefschetz trace formula.
We note that $R$ and $\Fr_q$ are far from being unique. Let $(R_1,\Fr_{q,1})$ be another lift. 
Then there is a canonical isomorphism $\psi$ between the cohomology groups associated to $\Omega_R^\bullet$ 
and the cohomology groups associated to $\Omega_{R_1}^\bullet$ such  that $\psi\circ \Fr_q=\Fr_{q,1}\circ \psi$.

\begin{example} Let $V$ be the curve given by $F=y^2z-x^3-axz^2-bz^3$ for $a,b\in \FF_q$ such 
that $\gcd(6,q)=1$ and $4a^3+27b^2\neq 0$. Let $X=\PP^2\setminus V$. Then $H^0_{MW}(X,\QQ_q)\cong \QQ_q$ and $H^2_{MW}(X,\QQ_q)$ is generated by
\[ \frac{xyz}{F}\Omega \mbox{ and } \frac{(xyz)^2}{F^2}\Omega\]
with $\Omega =dy\wedge dz-dx\wedge dz +dx\wedge dy$. All other cohomology groups vanish.

As a lift of Frobenius one can take $[x,y,z]\mapsto [x^p,y^p,z^p]$.
\end{example}
\begin{remark}
Suppose $E\subset \PP^2_{\CC}$ is a complex elliptic curve given by the zero-set of a cubic $F$ as in the above example. Then one can show that the usual de Rham cohomology group $H^2(\PP^2\setminus E,\CC)$ is generated by
\[ \omega_1:=\frac{XYZ}{F}\Omega \mbox{ and } \omega_2:=\frac{(XYZ)^2}{F^2}\Omega.\]
A famous theorem of Griffiths~\cite{GrRa} applied to this case yields that $\omega_2$ generates $F^2H^2(\PP^2\setminus E,\CC)$  and $F^1H^2(\PP^2\setminus E,\CC)/F^2H^2(\PP^2\setminus E,\CC)$ is generated by $\omega_1$, where $F^{\bullet}$ is the Hodge filtration.
\end{remark}

\begin{example} Let $g$ be a positive integer. Consider the algebra
\[ \overline{R}:=\FF_q[x,y,z]/(-y^2+f_{2g+1},yz-1)\]
with $f_{2g+1}\in \FF_q[x]$ a polynomial of degree $2g+1$, with only simple roots. 
Then $\Spec\overline{R}$ is a hyperelliptic curve of genus $g$ minus the $2g+2$ Weierstrass points. 
Lift $\overline{0},\overline{1},\dots,\overline{p-1}$ to $0,1, \dots ,p-1$. This choice induces a 
unique lift $\sigma: \QQ_p\rightarrow \QQ_p$ of Frobenius.

A basis for $H^1_{MW}(X,\QQ_q)$ is
\[ x^i dx \mbox{ for } i=0,\dots,2g, \mbox{ and } x^i\frac{dx}{y} \mbox{ for } i=0,\dots 2g-1.\]
There is a unique lift $\Fr_q$ of Frobenius extending $\sigma$, such 
that $\Fr_q(x)=x^q$, $\Fr_q(y)\equiv y^q\bmod \pi$ and $\Fr_q(y)^2=\Fr_q(f_{2g+1})$. Finally $\Fr_q(z)$ is 
defined by $\Fr_q(z)\Fr_q(y)=1$.
This description of the Monsky-Washnitzer cohomology enabled Kedlaya \cite{Ked} to present an 
efficient point counting algorithm on hyperelliptic curves over fields of small characteristic. 
We would like to point out that $\Fr_q(y)$ and $\Fr_q(z)$ are overconvergent power series which do not 
converge to a rational function. This is in contrast to the previous example.
\end{example}

The fact that we are working with a cohomology theory using differential forms enables us to use more 
analytic techniques. Consider the family $V_t$ of elliptic curves $X^3+Y^3+Z^3+3tXYZ$ where $t$ is a 
parameter in  the $p$-adic unit disc $\Delta$. Then $H^2_{MW}(X_{t_0},\QQ_q)$ is two-dimensional for a 
general $t_0\in \Delta$. One of the two eigenvalues of $p^2\Fr_p^{-1}$ has $p$-adic valuation $0$. 
This eigenvalue is called the unit root. 
Dwork showed that the unit root, as a function of $t$, satisfies the 
$p$-adic hypergeometric differential 
equation with parameters $\frac{1}{2},\frac{1}{2}; 1$. Specialising to values $t_0$ on the boundary of the 
unit disc gives the characteristic polynomial of Frobenius on 
$H^2_{MW}(X_{\overline{t_0}},\QQ_q)$ with $\overline{t_0}\equiv t \bmod \pi$.
N. Katz has generalised this idea: on the $p$-adic cohomology 
groups $H^{n+1}_{MW}(X_t)$ one can introduce a Gauss-Manin connection, much in the same way as one 
usually does in the complex case.

\begin{theorem}[{Katz, \cite{Katz}}] \label{KatzThm}Let $X_t$ be a family of hypersurface complements in  $\PP^n$ and let $(H^{n}_{MW}(X_t,\QQ_q),\nabla_{GM})$ be the local system of cohomology groups. 
Let $A(t)$ be a solution of the Picard-Fuchs equation (i.e., the differential equation 
associated to $\nabla_{GM}$.) Fix a basis $e_{1,t},\dots,e_{m,t}$ for this local system and let $\Fr(t)$ be 
the matrix of 
the Frobenius action on $H^{n}_{MW}(X_t,\QQ_q)$ with respect to $e_{1,t},\dots,e_{m,t}$ . 
Let $t_0\in \QQ_q$ such that $\mid t_0\mid_q=1$.  Then
\[ \Fr(t_0)=A(t_0^q)^{-1} \Fr(0) A(t_0)\]
provided that the reduction $\overline{X}_{\overline{t_0}}$ is smooth.
\end{theorem}

\subsection{Mirror Symmetry and Zeta Functions} \label{sub_zetamirror}
The classical Picard-Fuchs equations of Calabi-Yau threefolds play an important role in Mirror Symmetry.
Candelas, de la Ossa and Rodriguez-Villegas \cite{COR1}, \cite{COR2} considered the $p$-adic 
Picard-Fuchs equation (cf. Theorem~\ref{KatzThm}) of the following one-parameter family of
quintics
\begin{equation}\label{eqn_M}
\cM: \quad Q(\psi)= \left\{ \sum_{i=0}^{4} x_i^5 - 5 \psi x_0x_1x_2x_3x_4=0 \right\}.
\end{equation}
For general $\psi$, i.e., $\psi^5 \neq 0,1,\infty$, this is a smooth quintic. 
They calculated explicit expressions for the number of points on $Q(\psi)$, using tedious computations 
with so-called Dwork characters. Then they observed that these expression were solutions to the 
Picard-Fuchs equations. The above mentioned Theorem~\ref{KatzThm} gives a more mathematical 
explanation for this phenomenon.

Using these calculations,  
Candelas et al.  determined 
the structure of the zeta-function, more precisely, they showed that for any prime $p\neq 5$ and any $\psi\in \FF_{p} \setminus\{0,1\}$
\begin{equation}\label{eqn_zetaM}
Z_{\cM}(t,\psi)=Z_{Q(\psi)}(t)=
\frac{R_0(t,\psi) R_{A}(p^{\rho}t^{\rho},\psi)^{\frac{20}{\rho}} 
R_{B}(p^{\rho}t^{\rho},\psi)^{\frac{30}{\rho}} }
{(1-t)(1-pt)(1-p^2t)(1-p^3t)}
\end{equation}
where $\rho=1,2$ or $4$ is the least integer such that $p^{\rho} \equiv 1 \bmod 5$.
The $R$'s are quartic polynomials in $t$.

We would like to point out that this factorisation can be obtained in a more direct way. As far as we know this strategy is not contained in the literature.
Every member of the family $\cM$ has many automorphisms. We distinguish three types. 
First of all, each $\sigma\in S_5$ induces a  permutation of coordinates $x_i\mapsto x_{\sigma(i)}$. 
Secondly,   let $\zeta$ be a primitive $5$-th 
root of unity. Then the group
\begin{equation}\label{eqn_groupaction}
G= \left\{ \diag(\zeta^{n_0},\zeta^{n_1},\zeta^{n_2},\zeta^{n_3},\zeta^{n_4}); 
\sum n_i \equiv 0 \bmod 5 \right\} 
\end{equation}
acts on the family $\cM$. (Note that this action is defined over $\FF_{p^\rho}$.) 
Finally, the third type is the Frobenius automorphism $\Fr_{p^{\rho}}$. 
Note that $\Fr_{p^\rho}$ commutes with all elements in $G$ as well as all elements in $S_5$. It is easy to write down an explicit basis for $H^4_{MW}(\PP^4\setminus Q(\psi),\QQ_p)$ and to calculate the action of the automorphisms of the first and second type on this basis.

Since two commuting automorphisms respect each other's eigenspaces, we can use the automorphisms coming 
from the $S_5$-action to identify all the eigenspaces of Frobenius on  $H^4_{MW}(\PP^4\setminus Q(\psi,\QQ_p))$.
Then we can use the automorphisms of $G$ to prove that many of the eigenvalues of Frobenius coincide. This gives a factorisation of $Z_{\PP^4\setminus Q(\psi)}(t)$. Since $Z_{\PP^4\setminus Q(\psi)}(t)Z_{\cM}(t,\psi)=Z_{\PP^4}(t)$ this factorisation  can be translated in a factorisation of $Z_{\cM}(t,\psi)$ as in (\ref{eqn_zetaM}).

The $R$'s are quartic polynomials, satisfying the functional equation
$$
R\left(\frac{1}{p^3t},\psi\right)=\frac{1}{p^6t^4}R(t,\psi).
$$
In particular
$$
R_0(t,\psi)= 1 + a_1(\psi)t + b_1(\psi)pt^2 + a_1(\psi)p^3t^3 + p^6t^4
$$
where $a_1(\psi), b_1(\psi)$ depend on the parameter $\psi$. The two factors $R_A$ and $R_B$ are 
closely related to genus $4$ curves $A=A(\psi)$ and $B=B(\psi)$, namely
$$
Z_A(u,\psi)= \frac{R_A(u,\psi)^2}{(1-u)(1-pu)},
$$ 
with an analogous expression for $R_B$. The curves $A$ and $B$ are the nonsingular models of 
\[ y^5=x^2(1-x)^\beta(1-x/\psi^5)^{5-\beta}\]
where $\beta=3$ for $A$, and $\beta=4$ for $B$. The geometric role of $A$ and $B$ does not seem to be well understood.

There is no good notion of a mirror family over non-algebraically closed fields in positive 
characteristic. However, in this case we can mimic the mirror construction since the mirror 
family of (\ref{eqn_M}), considered as a family of complex varieties, is given by a quotient 
construction using $G$:
for a general element $Q(\psi)$ the quotient $Q(\psi)/G$ has nodes
and resolving these nodes gives Calabi-Yau threefolds $Q'(\psi)$, which form the mirror 
family $\cW$ of $\cM$. 

In the complex case
one has an equivalent construction, namely as family of
hypersurfaces in the toric variety $\PP_{\Delta}$ where $\Delta$ is the simplex in $\RR^4$ given by
$\{e_1, \ldots, e_4, -(e_1 + \cdots e_4)\}$. The mirror family is then 
given by crepant resolutions of the closure in $\PP_{\Delta}$ of the affine family in the torus $(\CC^*)^4$ 
defined by the equation
$$
x_1 + \cdots + x_n + \frac{1}{x_1 \cdots x_4} + \psi =0.
$$
It is known among experts how one can generalise this construction to $p$-adic toric varieties.

Candelas et al. computed the zeta function of the mirror family for general $\psi$ as
\begin{equation}\label{eqn_zetaW}
Z_{\cW}(t,\psi)=Z_{Q'(\psi)}(t)= \frac{R_0(t,\psi)}{(1-t)(1-pt)^{101}(1-p^2t)^{101}(1-p^3t)}.
\end{equation}

Using the explicit differential forms provided by Mon\-sky-Wash\-nit\-zer cohomology 
one can easily find an explicit basis for $H^4_{MW}(\PP^4\setminus Q(\psi))^G$. Along the same lines one 
shows that $R_0(t,\psi)$ is the characteristic polynomial of Frobenius on $H^4_{MW}(\PP^4\setminus Q(\psi))^G$. This explains why $R_0(t,\psi)$ is a factor of the numerator of the Zeta function of the mirror threefold.

Comparing equations (\ref{eqn_zetaM}) and  (\ref{eqn_zetaW}) shows that there is a close relationship
between the zeta function of the family $\cM$ and its mirror $\cW$. In particular, one obtains the
congruence
\begin{equation}\label{eqn_congruence}
Z_{\cM}(t,\psi) \equiv \frac{1}{Z_{\cW}(t,\psi)}\equiv (1-pt)^{100}(1-p^2t)^{100} \bmod 5^2.
\end{equation}

In her thesis Kadir \cite{Ka} has investigated the case of double octics, where she was able to prove
analogous results.

The question of studying the zeta functions of mirror families was also taken up by Wan \cite{Wa}
and Fu and Wan \cite{FW}. Wan considered the family of hypersurfaces 
$$
X_{\psi}: x_0^{n+1} + \cdots + x_n^{n+1} + \psi x_0 \cdots x_n= 0.    
$$
The construction of the `mirror family' $Y_{\psi}$ is then completely analogous 
to the quintic case. He calls a pair $(X_{\psi},Y_{\psi})$ with the same parameter $\psi$ a 
\emph {strong} mirror pair and proves the following generalisation of (\ref{eqn_congruence}) for such pairs:
\begin{equation}\label{eqn_congruence2}
\# X_{\psi}(\FF_{q^k}) \equiv \# Y_{\psi}(\FF_{q^k}) \bmod q^k.
\end{equation}
Here $q=p^r$ is a prime power.
He conjectures that such congruences should hold in greater generality. Such a generalisation for 
certain pairs $(X,X/G)$ was then indeed obtained by Fu and Wan in \cite{FW}. As a consequence of their
results they prove
\begin{theorem}
Let $X$ be a smooth, geometrically connected variety $X$ over the finite field $\FF_q$, $q=p^r$.
Assume that $H^i(X,\cO_X)=0$ for $i>0$. Then for all integers $k$, the congruence
$$
\# X(\FF_{q^k}) \equiv 1 \bmod q^k
$$
holds.  Moreover, let $G$ be a finite group of $\FF_q$ automorphisms acting on $X$. Then 
$$
\# (X/G)(\FF_{q^k}) \equiv \# X(\FF_{q^k}) \equiv 1 \bmod q^k.
$$   
\end{theorem} 
Their proof uses crystalline cohomology and the Mazur-Ogus theorem. Esnault \cite{Es} has pointed out that 
these results can also be obtained using de Rham-Witt cohomology and rigid cohomology.

\section{Open Questions}\label{Sec_open}

We conclude this paper with a brief discussion of some open problems. 
These will mainly concern Calabi-Yau threefolds. We start with some questions related to modular forms.

If $f$ is a newform of weight two with field of coefficients $F$, then a classical construction 
associates to $f$ an abelian variety $A$ over $\QQ$ of dimension $[F:\QQ]$. 
In particular, if $f$ has rational coefficients, then this gives an elliptic curve $E$. 
For higher weight, however, no such construction is available. 
As remarked in Remark~\ref{rem_repgeo} one can always find some symmetric product of some universal curve for which
part of the cohomology corresponds to $f$. One might ask whether any such form can 
be realised geometrically by a Calabi-Yau variety:

\begin{question}\label{ques_newforms}
Let $f$ be a newform of weight $k>2$ with rational Fourier coefficients. 
Is there a Calabi-Yau variety $X$ over $\QQ$ such that $L(f,s)$ occurs in the $L$-series of $X$?
\end{question}
This question was formulated independently  by B. Mazur and D. van Straten.

Consider two non-isomorphic modular varieties defined over $\QQ$ such that they are isomorphic over $\overline{\QQ}$. In practice, it turns out that very often the associated modular forms differ by a twist.
E.g., this can be seen in the case of 
elliptic fibrations in Weierstrass form or with double octics, where twists over quadratic extensions
can easily be realised. 
For this reason the above question should be asked modulo twisting.
By Theorem~\ref{theo_CMforms} (and assuming GRH) the number of available forms of given odd weight is finite 
up to twisting. In the case of $\Kthree$ surfaces Shafarevich's theorem
\ref{theo_Shafarevich} also tells us that the number of available surfaces is finite, hence in the weight three case
this reduces to a question of matching two finite sets. We have already mentioned that, so far, $24$ of the 
$64$ forms have been realised geometrically (see Subsection~\ref{sub_examplesK3} or \cite[Section 15]{Sch3}).

On the other hand, if the weight $k$ is even, then there are infinitely many newforms with rational coefficients, 
even up to twisting. An answer to Question~\ref{ques_newforms} would therefore require a general 
construction, similar to the weight two case of elliptic curves. In the case of weight four, we have seen 
a  number of examples of modular 
Calabi-Yau threefolds in Section~\ref{sub_rigidCY}. In total, Calabi-Yau threefolds for around $80$ newforms 
have been found so far (modulo twisting). We refer the reader to Meyer's book  \cite{CM}, which 
contains the most systematic collection of examples up to date. 
These include non-rigid and nodal varieties as well as some 
threefolds where the corresponding newform is only conjectural. 
{}Nevertheless, there are no known Calabi-Yau partners for  the newforms of level $7$ 
and $13$. 
In general, it seems to be hard to find modular varieties with 
big bad primes. It would be interesting to find constructions which yield such examples.

Another question concerns the level of the corresponding newform:

\begin{question}
{Let $X$ be a  modular variety over $\QQ$. Can one predict the level of the associated newform?}
\end{question}

In particular, this question concerns the bad primes of $X$. 
Once more, the difference with the special case of elliptic curves has to be emphasised. 
For an elliptic curve $E$, say over $\QQ$ and hence modular, the theory provides a globally minimal model. 
This can be used to define the conductor of the curve, which then coincides with the level of the corresponding 
newform.

For a general variety $X$ over $\QQ$, no concept of arithmetic minimality is known. 
Nevertheless, we can define a conductor for any compatible system of \'etale Galois representations. 
Assuming modularity, this conductor coincides with the level $N$ of the associated newform. 
This is a consequence of the compatibility with the local Langlands correspondence. 

In Section~\ref{sub-Wiles}, we gave the classical bounds for the exponents $e_p$ of the conductor, if $E$ is an 
elliptic curve.
The same bounds also hold in many other cases. 
For a smooth projective variety $X$ over $\QQ$ of odd dimension $n$, these were proved by 
Serre in \cite{Se2}
for the Galois representation of $H^n_{\etal}(X,\QQ_l)$ 
provided this is two-dimensional and the Hodge decomposition is
\[
H^n(X_{\CC}, \CC)\cong H^{n,0}(X_\CC)\oplus H^{0,n}(X_\CC),
\]
and assuming that the representation is modular. This was reproved by Dieulefait in \cite{Di1}.
Dieulefait also remarks that his 
proof (combined with the potential modularity result of Taylor) is enough to 
conclude that the result still holds without the modularity assumption.
In particular, these bounds apply to rigid Calabi-Yau threefolds over $\QQ$. For newforms of odd weight and 
rational coefficients (or in general CM-forms), the same bounds have been obtained by Sch\"utt 
in \cite[Cor. 12.2]{Sch3}. In practice, these bounds can be used to precisely determine the  
corresponding newform, 
if modularity is known (e.g., by Theorem~\ref{theo:DM}). One then writes down all newforms with level 
composed of the bad primes and checks some coefficients after counting points. 
For more details we refer to \cite{Di1}, where this was discussed for Calabi-Yau varieties.
At the moment, this method only works if there are at most two very small bad primes 
or if the newform has CM.

This takes us back to singular $\Kthree$ surfaces and rigid 
Calabi-Yau threefolds over $\QQ$. If a bad prime $p$ divides the level of the newform 
associated to $X$, then examples show the following: if the reduction $X_p$ has a single $A_1$-singularity 
 then $e_p=1$. Otherwise, $e_p= 2$ (for $p>3$) or  $e_p\geq 2$, if $p=2$ or $3$). In view of this, we would like to 
ask the following question:

\begin{question}
Let $X$ be a modular variety over $\QQ$ with corresponding newform $f$. 
Under which conditions do  the bad Euler factors of $L(X,s)$ and $L(f,s)$ coincide?
\end{question}

Since the conductors have to be equal, one can hope for an unconditional answer. 
To our knowledge, this has not yet been investigated.

We shall now mention one further tool which can sometimes be used to prove modularity. 
This consists in exhibiting a birational map (over $\QQ$) between a modular variety and one which is 
expected to be 
modular. In \cite{HSvGvS}, two such varieties were called \emph{relatives}. 
This paper gives a number of examples for such relatives.
For example, birational maps over $\QQ$  between the following varieties were given: the Barth-Nieto quintic, 
the self-fibre product of the universal elliptic curve $\mathcal{E}_1(6)$ and Verrill's threefold associated to 
the root lattice $A_3$. 
Recently, Cynk and Meyer published a study of rigid Calabi-Yau threefolds of level $8$ where the question of
relatives was
treated systematically \cite{CM2}.

This issue can also be discussed for the self-fibre products from Section~\ref{subsub_fibre}. Here, each permutation of four cusps leads to isogenous 
elliptic fibrations. Then it is a matter of the field of definition of the 
isogeny which newform occurs. If the isogeny is not defined over $\QQ$, then 
the newform might be twisted. For instance, the modular elliptic surfaces 
associated to $\Gamma_1(4)\cap\Gamma(2)$ and $\Gamma_0(8)\cap\Gamma_1(4)$ are 
isogenous. If one works with the equations from \cite{B}, this isogeny is 
only defined over $\QQ(i)$. The newforms of level $8$ and $16$ are 
twists by the corresponding quadratic character.

In general, this problem is addressed by the \emph{Tate conjecture}. 
Let $\rho_1, \rho_2$ be isomorphic two-dimensional $l$-adic Galois representations coming from geometry. 
In particular, these Galois representations occur in the \'etale cohomology of some smooth projective 
varieties $X_1, X_2$ over $\QQ$. 
Then  the Tate conjecture predicts a correspondence, i.e., an algebraic cycle in the 
product $X_1\times X_2$, which is defined over $\QQ$ and which 
induces an isomorphism between $\rho_1$ and $\rho_2$.

For details on known correspondences, the reader is referred to the book of Meyer \cite{Me}. 
Here we shall only mention one particular example. In \cite{vGN}, Van Geemen and Nygaard study three 
Siegel modular threefolds over $\QQ$. One of them is the first non-rigid Calabi-Yau threefold for which 
modularity has been established. Denoting this variety by $X$, Van Geemen and 
Nygaard show that $h^{1,2}(X)=1$.
An automorphism of $X$ is used to split the four-dimensional Galois representation of 
$H^3_{\etal}(X,\QQ_l)$ into a sum of two-dimensional representations. 
Then the authors deduce the equality
\[
L(X,s)=L(f_4,s)\,L(f_2,s-1)
\]
with newforms of level $32$ and weight four and two respectively. 
Here, $f_2$ comes from the Gr\"ossencharacter $\psi$ of 
$\QQ(i)$ with conductor $(2+2i)$ and $\infty$-type $1$. 
The other newform, $f_4$ is derived from $\psi^3$. 
On the other hand, $\psi$ is related to the elliptic curve with CM by $\ZZ[i]$:
\[
E:\;\; y^2=x^3-x \;\;\text{ as }\;\;\; L(E,s)=L(\psi,s).
\]
As a consequence, 
\[
L(X,s)=L(\Sym^3 H^1_{\etal}(E,\QQ_l), s),
\]
and hence the Tate conjecture predicts the existence of a correspondence between $X$ and $E\times E\times E$ 
which induces an isomorphism between the two representations corresponding to $\psi^3$. 
This is established by means of a curve of genus $5$ over $\QQ$ which admits rational maps over $\QQ$ to 
both varieties, $X$ and $E\times E\times E$.

However, there are modular varieties $X$ and $Y$ with the same associated form, for which there is no known correspondence.

\begin{question} Let $X$ and $Y$ be modular varieties which have an  associated form in common, find a correspondence between them.
\end{question}

We would like to conclude this section with a question related to arithmetic mirror symmetry.
In Section~\ref{sub_zetamirror} we discussed the zeta function of the family $\cM$ of quintics
in $\PP^4$ (See~(\ref{eqn_zetaM})). This expression involves genus $4$  curves $A$ and
$B$ whose geometric meaning has, so far, remained a mystery.   

\begin{question}
What is the geometric meaning of the genus $4$ curves which appear in the motive of the
one-parameter quintic family $\cM$?
\end{question}

\end{document}